\documentclass[11pt]{amsart}

\usepackage{amsmath,amssymb,amsthm,mathtools}
\usepackage{enumitem}
\usepackage{microtype}
\usepackage[hidelinks]{hyperref}
\usepackage{fullpage}

\newtheorem{theorem}{Theorem}[section]
\newtheorem{corollary}[theorem]{Corollary}
\newtheorem{lemma}[theorem]{Lemma}
\newtheorem{proposition}[theorem]{Proposition}
\newtheorem{conjecture}[theorem]{Conjecture}
\theoremstyle{remark}
\newtheorem{remark}[theorem]{Remark}

\newcommand{\PP}{\mathbb P}
\newcommand{\kk}{\Bbbk}
\newcommand{\NN}{\mathbb N}
\newcommand{\ZZ}{\mathbb Z}
\newcommand{\Spec}{\operatorname{Spec}}

\newcommand{\Hilb}{\operatorname{Hilb}}
\newcommand{\vol}{\operatorname{vol}}
\newcommand{\eul}[2]{\genfrac{\langle}{\rangle}{0pt}{}{#1}{#2}}

\title{Frobenius Dilation of Hilbert Functions, Demailly's Conjecture, and Varieties of Powers}
\author{Trung Hoa Dinh}
\address{Department of Mathematics and Statistics, Troy University, Troy, Alabama 36082, USA}
\email{tdinh@troy.edu}

\author{T\`ai Huy H\`a}
\address{Tulane University, Department of Mathematics, 6823 St. Charles Avenue, New Orleans, Louisiana 70118, USA}
\email{tha@tulane.edu}

\subjclass[2020]{14N07, 13A35, 14N20, 13D02}
\keywords{Frobenius powers, Hilbert functions, Waring problem, powers of forms, generic $k$-rank, secant varieties, Seshadri constants, apolarity, fat points, symbolic powers, Waldschmidt constants, Demailly conjecture}

\begin{document}

\begin{abstract}
We prove a Frobenius-dilation inequality for Hilbert functions of ideals generated by powers of fixed homogeneous forms in characteristic zero. Via Emsalem--Iarrobino inverse system duality, it yields quantitative fat-point postulation estimates and a new proof of the recent theorem of H\`a--Sivakumar resolving Demailly's conjecture for arbitrary finite point sets. A fixed-characteristic version controls tangent spaces to varieties of powers and, together with plane interpolation and secant-theoretic results, proves complete nondefectivity of $V^k_{2,d}$ for all $k\ge3$ and $d\ge2$; in particular, it settles the ternary Waring problem for generic $k$-rank conjectured by Lundqvist--Oneto--Reznick--Shapiro (suggested by Ottaviani). We also show that, for a prime power exponent $q$, any failure of Nicklasson's conjecture in three variables is confined to $2q-2$ consecutive degrees.
\end{abstract}

\maketitle


\section{Introduction}\label{sec:intro}

Throughout the paper, $\kk$ is an algebraically closed field of characteristic zero and $\NN=\{1,2,\ldots\}$. Let
$
S=\kk[y_1,\ldots,y_N],
$
let $g_1,\ldots,g_s\in S$ be homogeneous forms, and set
$$
J=(g_1,\ldots,g_s),
\quad
J_q=(g_1^q,\ldots,g_s^q).
$$
The ideal $J_q$ is generated by the $q$th powers of the chosen forms. In characteristic zero it generally differs from the ordinary power $J^q$ and depends on the chosen generating sequence. Ideals of this type also occur in work of Ma--Smirnov \cite{MaSmirnov} on refinements of Lech's inequality.

When $q=p^e$, the situation in characteristic $p$ is rigid. If $\kappa$ is a perfect field of characteristic $p$, $T=\kappa[y_1,\ldots,y_N]$, and $I\subseteq T$ is homogeneous, then Frobenius gives
$$
\Hilb_{T/I^{[q]}}(z)
=
\Hilb_{T/I}(z^q)(1+z+\cdots+z^{q-1})^N.
$$
There is no Frobenius endomorphism in characteristic zero relating $J$ and $J_q$. The starting point of this paper is that a coefficientwise analog of the characteristic-$p$ identity survives in characteristic zero.

\medskip
\noindent\textbf{Frobenius dilation of Hilbert functions.}
Our starting observation is Theorem \ref{thm:HF-dilation}, which shows that there is a finite set $\Sigma$ of rational primes such that, for every $p\notin\Sigma$ and every $q=p^e$,
\begin{equation}\label{eq:intro-HF}
\Hilb_{S/J_q}(z)
\preceq
\Hilb_{S/J}(z^q)(1+z+\cdots+z^{q-1})^N.
\end{equation}
A finite monic Gr\"obner basis is used to preserve the entire Hilbert function after reduction. In characteristic $p$ the right-hand side of \eqref{eq:intro-HF} is the exact Hilbert series of the corresponding Frobenius power, and semicontinuity of matrix rank gives the coefficientwise inequality after lifting back to characteristic zero. An important point is that one finite set of good primes works simultaneously in every degree and for every power of a fixed good prime.

Corollary \ref{thm:frob-amp} is the single-degree consequence that will be used repeatedly. If $(g_1,\ldots,g_s)_A=S_A$, then for every good prime power $q$,
$$
(g_1^q,\ldots,g_s^q)_{qA+(N-1)(q-1)}
=
S_{qA+(N-1)(q-1)}.
$$
The paper then develops this dilation principle in polynomial interpolation, in secant problems for varieties of powers, and in the Hilbert-series problem for powers of general forms.

\medskip
\noindent\textbf{Polynomial interpolation and Demailly's conjecture.}
Our first application concerns symbolic powers of ideals of points. Let $Z\subset\PP_\kk^n$ be a nonempty finite reduced set of points with defining deal $I=I(Z)$. For a homogeneous ideal $K$, let $\alpha(K)$ denote its initial degree, and let
$$
\widehat\alpha(K)
=
\lim_{t\to\infty}\frac{\alpha(K^{(t)})}{t}
$$
be its Waldschmidt constant. Demailly \cite{Demailly} conjectured in 1982 that, for every $m\ge1$,
\begin{equation}\label{eq:intro-Demailly}
\widehat\alpha(I)
\ge
\frac{\alpha(I^{(m)})+n-1}{m+n-1}.
\end{equation}
For $m=1$ this is Chudnovsky's conjecture \cite{Chudnovsky}. The conjecture for arbitrary finite point sets was resolved recently by H\`a--Sivakumar \cite[Theorem~1.2]{HaSivakumar}. We obtain a different proof from Theorem \ref{thm:HF-dilation}, and the proof arises from a more quantitative comparison of fat-point postulation functions.

The bridge is the Emsalem--Iarrobino inverse system correspondence \cite[Theorem~I]{EI}. For distinct points $P_1,\ldots,P_s\in\PP^n_\kk$ and a multiplicity vector $\mathbf u=(u_1,\ldots,u_s)$, write
$$
I(\mathbf u)=\bigcap_{i=1}^s I(P_i)^{u_i},
\quad
h_{\mathbf u}(D)=\dim_\kk I(\mathbf u)_D,
$$
with the conventions of Section \ref{sec:fatpoints}. Theorem \ref{thm:fatpoint-HF} translates Theorem \ref{thm:HF-dilation} into a coefficientwise inequality comparing the postulation of a Frobenius-dilated multiplicity vector with finitely many neighboring postulation values of the original vector. Its vanishing specialization is Theorem \ref{thm:fatpoint-amp}. If
$
a=\alpha\left(\bigcap_i I(P_i)^{u_i}\right),
$
then, for every good prime power $q$ and every positive vector $\mathbf v$ satisfying
$$
v_i\ge qu_i+(n-1)(q-1),
$$
one has
\begin{equation}\label{eq:intro-fatpoint-amp}
\alpha(I(\mathbf v))\ge qa+(n-1)(q-1).
\end{equation}
For uniform multiplicities, \eqref{eq:intro-fatpoint-amp} gives Corollary \ref{cor:points}, namely Demailly's inequality \eqref{eq:intro-Demailly}, after dividing by the amplified multiplicity and letting $q=p^e$ tend to infinity.

The same argument is not restricted to uniform symbolic powers. Theorem \ref{thm:main} treats positive subadditive multiplicity sequences $\mu_i(t)$ with asymptotic slopes $\lambda_i$ and, for
$$
F_t=\bigcap_i I(P_i)^{\mu_i(t)},
$$
gives
$$
\widehat\alpha(\mathcal F)
\ge
\min_i
\frac{\lambda_i\bigl(\alpha(F_m)+n-1\bigr)}{\mu_i(m)+n-1}, \text{ where } \lambda_i = \lim_{t \rightarrow \infty} \frac{\mu_i(t)}{t}.
$$

The postulation inequality also has an asymptotic geometric consequence. In Theorem \ref{thm:Eulerian-volume} we give bounds for the volume of certain divisor on the blow-up of $\PP^n$ at the points by a finite Eulerian-weighted combination of neighboring postulation values. Thus the same Hilbert-function inequality that yields Demailly's conjecture also controls a natural asymptotic positivity invariant on blow-ups.

\medskip
\noindent\textbf{Varieties of powers and the ternary Waring problem.}
A second application concerns Waring-type decompositions in which the summands are powers of forms of a prescribed degree. Geometrically, these decompositions are governed by secant varieties of varieties of powers. For integers $k,n,d\ge2$, set
$$
V^k_{n-1,d}
=
\{[G^k]:G\in\kk[x_1,\ldots,x_n]_d\}
\subseteq
\PP(\kk[x_1,\ldots,x_n]_{kd}).
$$
For a projective variety $V\subseteq\PP^M$, let $\sigma_s(V)$ denote its $s$th secant variety, the Zariski closure of the union of the linear spans of $s$ points of $V$. Its expected dimension is
$$
\min\{s(\dim V+1)-1,M\}.
$$
We call $V$ \emph{completely nondefective} if every secant variety $\sigma_s(V)$ has the expected dimension.

The $k$-rank of a form $F\in\kk[x_1,\ldots,x_n]_{kd}$ is the least integer $s$ for which
$$
F=G_1^k+\cdots+G_s^k,
$$
with $G_i\in\kk[x_1,\ldots,x_n]_d$. The value for a general form is denoted by $\operatorname{rk}^{\circ}_k(n,kd)$; equivalently, it is the least $s$ for which $\sigma_s(V^k_{n-1,d})$ fills its ambient projective space. Lundqvist--Oneto--Reznick--Shapiro \cite[Conjecture~1.2]{LORS}, following a suggestion of Ottaviani, conjectured that, for $k\ge3$,
$$
\operatorname{rk}^{\circ}_k(n,kd)
=
\min\left\{s\ge1:
 s\binom{d+n-1}{n-1}\ge \binom{kd+n-1}{n-1}\right\}.
$$
The conjecture is known for binary forms, but remains open in general in higher dimension.

Theorem \ref{thm:syzygy-gap} gives the first algebraic input for the secant problem: if $\delta(g_1,\ldots,g_s)$ denotes the least coefficient degree of a nontrivial syzygy among forms of the same degree, then
$$
\delta(g_1^q,\ldots,g_s^q)
\ge q\,\delta(g_1,\ldots,g_s)
$$
for every good prime power $q$. The fixed-characteristic counterpart is Theorem \ref{thm:secant-tangent}, which gives a numerical criterion for injectivity of the tangent map to a secant variety of powers. Corollary \ref{cor:frobenius-complete} then exhibits a uniform large-degree consequence: if $k-1$ is a prime power and $d\ge3k^3-6k^2+4$, then $V^k_{2,d}$ is completely nondefective. In particular, this gives complete nondefectivity for cubes when $d\ge31$ and for fourth powers when $d\ge100$.

The cube case can be completed in all degrees by reducing to considering linear systems on blow-ups of $\PP^2$ at at most eight general points, with two exceptional pairs treated separately. The resulting Theorem \ref{thm:ternary-cubes-all} proves complete nondefectivity of $V^3_{2,d}$ for every $d\ge1$, and Corollary \ref{cor:generic-ternary-cube-rank} gives the corresponding generic cube-rank formula.

For higher powers, Lemma \ref{lem:homogeneous-postulation} packages the homogeneous plane-interpolation results needed for multiplicities at most $42$, while Lemma \ref{lem:fourth-seshadri} supplies the multipoint Seshadri estimate used for fourth powers. Together with the Frobenius ranges above, the filling theorem of Fr\"oberg--Ottaviani--Shapiro \cite{FOS}, and the degeneration theorem of Casarotti--Postinghel \cite{CasarottiPostinghel}, these ingredients yield our main result on varieties of powers.

\medskip
\noindent\textbf{Theorem \ref{thm:ternary-k-rank}.}
\emph{For every $k\ge3$, $d\ge2$, and $s\ge1$,
$$
\dim\sigma_s(V^k_{2,d})
=
\min\left\{
 s\binom{d+2}{2}-1,
 \binom{kd+2}{2}-1
\right\}.
$$
Consequently,
$$
\operatorname{rk}^{\circ}_k(3,kd)
=
\left\lceil
\frac{\binom{kd+2}{2}}{\binom{d+2}{2}}
\right\rceil.
$$}
Thus the ternary case of the generic $k$-rank conjecture holds for every $k\ge3$; in fact, every ternary variety of powers $V^k_{2,d}$ with $k\ge3$ and $d\ge2$ is completely nondefective. Related recent work on powers of forms includes the identifiability results of Blomenhofer--Casarotti--Micha{\l}ek--Oneto \cite{BCMO} and the general nondefectivity bounds of Casarotti--Postinghel \cite{CasarottiPostinghel}.

\medskip
\noindent\textbf{Powers of general forms and Nicklasson's conjecture.}
The tangent-space arguments above require control only in a finite range of degrees. Nicklasson's conjecture \cite{Nicklasson} asks for the stronger assertion that the entire Hilbert series of $(G_1^m,\ldots,G_s^m)$ agrees with that of $s$ general forms of degree $dm$.

Theorem \ref{thm:Nenashev-range} combines Theorem \ref{thm:syzygy-gap} with Nenashev's maximal-rank theorem to give explicit degree ranges in which powers of general forms already have the generic Hilbert function. In three variables, Proposition \ref{prop:Nicklasson-two-maps} reduces Nicklasson's conjecture for fixed $(s,d,m)$ to two adjacent rank conditions, namely the last injective multiplication map and the first surjective one. Proposition \ref{prop:Nicklasson-four} proves the conjecture whenever $s\le4$, using the strong Lefschetz property for monomial complete intersections.

For prime-power exponents, Theorem \ref{thm:Nicklasson-three-vars} gives a more precise localization. Write
$$
\frac{(1-z^d)^s}{(1-z)^3}=\sum_{j\ge0} f_jz^j,
\quad
\tau(s,d)=\min\{j:f_j\le0\}.
$$
If $q=p^e$, then any possible failure of Nicklasson's conjecture can occur only in the $2q-2$ consecutive degrees
$$
q\tau(s,d)\le D\le q\tau(s,d)+2q-3.
$$

\medskip
\noindent\textbf{Further algebraic consequences.}
The coefficientwise inequality \eqref{eq:intro-HF} also yields two direct numerical consequences. Corollary \ref{cor:multiplicity} gives
$$
e(S/J_q)\le q^{\operatorname{ht}J}e(S/J)
$$
for every good prime power $q$, with equality when the chosen sequence is a homogeneous regular sequence of length $\operatorname{ht}J$. In the $\mathfrak m$-primary case, Corollary \ref{cor:Lech} gives
$$
\ell(S/J^q)+\ell(J^q/J_q)=\ell(S/J_q)\le q^N\ell(S/J),
$$
and the classical Lech inequality follows after passing to the limit.

\medskip
\noindent\textbf{Acknowledgment.} The second author is partially supported by a Simons Foundation grant.

\medskip
\noindent\textbf{Disclosure on the use of generative AI.} Generative AI was used during the development of this work as an assistive research tool. At an early stage, while we were searching for a second proof of Demailly's conjecture, it suggested examining our Hilbert-function results through the Emsalem--Iarrobino inverse system formula. This suggestion led us to investigate fat-point postulation and ultimately arrive at the new proof for Demailly's conjecture. It also suggested investigating whether our Hilbert-function methods might have consequences for Ottaviani's conjecture on generic $k$-rank, which led us to explore varieties of powers and their secant varieties. In addition, generative AI was used for language and structural editing and algebraic cross-checks. All mathematical statements, proofs, computations, and references were independently developed, examined and verified by the authors, who take full responsibility for the content of the manuscript.

\section{Frobenius dilation of Hilbert functions}\label{sec:frobenius}

Throughout this section let $S=\kk[y_1,\ldots,y_N]$. If $M$ is a finitely generated graded $S$-module, then write
$$
\Hilb_M(z)=\sum_{d\ge0}(\dim_{\kk}M_d)z^d.
$$
For two formal power series $A(z)=\sum a_dz^d$ and $B(z)=\sum b_dz^d$ with real coefficients, write $A(z)\preceq B(z)$ if $a_d\le b_d$ for every $d$. For $q\ge1$, set
$$
C_{N,q}(a)=[z^a](1+z+\cdots+z^{q-1})^N,
$$
with $C_{N,q}(a)=0$ when $a<0$ or $a>N(q-1)$.

\subsection{Reduction preserving Hilbert function}

To compare Hilbert functions after reduction, we need a single model that preserves all graded pieces at once. A finite monic Gr\"obner basis provides such a model: it replaces the infinitely many rank conditions, one in each degree, by finitely many polynomial identities that can be spread out over a finitely generated $\ZZ$-algebra.

\begin{lemma}\label{lem:HF-good-reduction}
Let $g_1,\ldots,g_s\in S$ be homogeneous and put $J=(g_1,\ldots,g_s)$. There exist a finitely generated integral $\ZZ$-subalgebra $B\subseteq\kk$ and a nonzero element $b\in B$ with the following property. For every field $\kappa$ and every homomorphism $B_b\to\kappa$, if $J_\kappa\subseteq\kappa[y_1,\ldots,y_N]$ denotes the ideal generated by the specialized forms, then
$$
\Hilb_{\kappa[y_1,\ldots,y_N]/J_\kappa}(z)=\Hilb_{S/J}(z).
$$
Moreover, after excluding finitely many rational primes, for every remaining prime $p$ there is a homomorphism $B_b\to\kappa_p$ to a finite field $\kappa_p$ of characteristic $p$.
\end{lemma}

\begin{proof}
Fix a monomial order and choose a finite homogeneous Gr\"obner basis $h_1,\ldots,h_t$ of $J$. After dividing by the leading coefficients, we may assume that the $h_j$ are monic. Choose a finitely generated integral $\ZZ$-subalgebra $B\subseteq\kk$ containing the coefficients of the $g_i$ and $h_j$, together with the coefficients occurring in expressions
$$
h_j=\sum_i a_{ji}g_i,
\quad
g_i=\sum_j c_{ij}h_j,
$$
and in fixed reductions of all $S$-polynomials of the $h_j$ to zero. Let $b\in B$ be the product of the finitely many nonzero denominators needed for these identities and pass to $B_b$.

Because the $h_j$ are monic, their leading monomials survive under every specialization. The chosen Buchberger reductions specialize as well, so the specialized $h_j$ form a Gr\"obner basis. The two displayed sets of identities show that they generate the same ideal as the specialized $g_i$. Thus the initial ideal is independent of the specialization and agrees with the initial ideal of $J$, so the graded Hilbert functions agree.

Since $B_b$ is a finitely generated integral $\ZZ$-algebra of characteristic zero, $\Spec B_b\to\Spec\ZZ$ is dominant. By Chevalley's theorem, its image contains a nonempty open subset of $\Spec\ZZ$, hence all but finitely many rational primes. For each such prime $p$, the fiber is a nonempty finite-type scheme over $\mathbb F_p$ and has a closed point. Its residue field $\kappa_p$ is finite by Zariski's lemma.
\end{proof}

\subsection{The characteristic-\texorpdfstring{$p$}{p} identity and characteristic-zero dilation}

In the following lemma, the exponent-class decomposition for Frobenius powers already appeared in \cite[Lemma~3.5]{HaSivakumar}.

\begin{lemma}\label{lem:HF-charp}
Let $\kappa$ be a perfect field of characteristic $p>0$, let $T=\kappa[y_1,\ldots,y_N]$, and let $I\subseteq T$ be homogeneous. For $q=p^e$ one has
\begin{equation}\label{eq:HF-charp}
\Hilb_{T/I^{[q]}}(z)
=
\Hilb_{T/I}(z^q)(1+z+\cdots+z^{q-1})^N.
\end{equation}
Equivalently, for every $D$,
\begin{equation}\label{eq:HF-charp-coeff}
\dim_\kappa(T/I^{[q]})_D
=
\sum_{m\ge0}C_{N,q}(D-qm)\dim_\kappa(T/I)_m.
\end{equation}
\end{lemma}

\begin{proof}
Because $\kappa$ is perfect, $T$ is free over $T^q$ with basis
$$
\{y^\beta:\beta\in\{0,\ldots,q-1\}^N\}.
$$
Equivalently, Frobenius is finite and flat. Tensoring $T/I$ along the $e$th Frobenius gives $T/I^{[q]}$. Keeping track of degrees yields a graded decomposition
$$
T/I^{[q]}
\cong
\bigoplus_{\beta\in\{0,\ldots,q-1\}^N}(T/I)^{(q)}(-|\beta|),
$$
where $(T/I)^{(q)}$ denotes the Frobenius twist of $T/I$, with all degrees multiplied by $q$. Since $\kappa$ is perfect, this twist has the same graded $\kappa$-dimensions as $T/I$. Taking Hilbert series and coefficients gives the two assertions.
\end{proof}

Recall that for $J=(g_1,\ldots,g_s)$ and $q\ge1$, write
$$
J_q=(g_1^q,\ldots,g_s^q).
$$
In characteristic zero this ideal depends on the chosen generating sequence. We now combine this reduction lemma with the exact characteristic-$p$ identity. 

\begin{theorem}\label{thm:HF-dilation}
Let $g_1,\ldots,g_s\in S$ be homogeneous and put $J=(g_1,\ldots,g_s)$. There is a finite set $\Sigma$ of rational primes such that, for every $p\notin\Sigma$ and every $q=p^e$,
\begin{equation}\label{eq:HF-dilation}
\Hilb_{S/J_q}(z)
\preceq
\Hilb_{S/J}(z^q)(1+z+\cdots+z^{q-1})^N.
\end{equation}
Equivalently, for every $D\ge0$,
\begin{equation}\label{eq:HF-dilation-coeff}
\dim_{\kk}(S/J_q)_D
\le
\sum_{m\ge0}C_{N,q}(D-qm)\dim_{\kk}(S/J)_m.
\end{equation}
\end{theorem}

\begin{proof}
Apply Lemma \ref{lem:HF-good-reduction}. After enlarging the finite set of bad primes, every remaining prime $p$ admits a finite-field specialization $B_b\to\kappa_p$ for which
$$
\Hilb_{T_p/J_p}(z)=\Hilb_{S/J}(z),
$$
where $T_p=\kappa_p[y_1,\ldots,y_N]$ and $J_p$ is generated by the specialized forms. Fix such a prime $p$ and let $q=p^e$. Since reduction commutes with ordinary $q$th powers and $q$ is a $p$th power, the specialization of $J_q$ is
$$
(J_q)_p=(\overline g_1^{\,q},\ldots,\overline g_s^{\,q})=J_p^{[q]}.
$$
Lemma \ref{lem:HF-charp} therefore gives
$$
\Hilb_{T_p/J_p^{[q]}}(z)
=
\Hilb_{S/J}(z^q)(1+z+\cdots+z^{q-1})^N.
$$

Fix $D\ge0$. The degree-$D$ piece of $J_q$ is the image of
$$
\Phi_{q,D}:\bigoplus_{i=1}^sS_{D-q\deg g_i}\longrightarrow S_D,
\quad
(h_i)_i\longmapsto\sum_i h_i g_i^q.
$$
With respect to monomial bases, $\Phi_{q,D}$ is represented by a matrix over $B_b$, whose specialization is the corresponding map over $\kappa_p$. Hence
$$
\operatorname{rank}_{\kk}\Phi_{q,D}
\ge
\operatorname{rank}_{\kappa_p}(\Phi_{q,D})_p.
$$
Equivalently,
$
\dim_{\kk}(J_q)_D
\ge
\dim_{\kappa_p}(J_p^{[q]})_D,
$
and therefore
$$
\begin{aligned}
\dim_{\kk}(S/J_q)_D
&=\dim_{\kk}S_D-\dim_{\kk}(J_q)_D\\
&\le \dim_{\kappa_p}(T_p)_D-\dim_{\kappa_p}(J_p^{[q]})_D\\
&=\dim_{\kappa_p}(T_p/J_p^{[q]})_D.
\end{aligned}
$$
Using \eqref{eq:HF-charp-coeff} and the equality of Hilbert functions gives \eqref{eq:HF-dilation-coeff}, hence \eqref{eq:HF-dilation}.
\end{proof}

\medskip

\noindent\textbf{Notation.} For the rest of the paper, a prime outside the finite exceptional set in Theorem \ref{thm:HF-dilation}, and any of its powers $q=p^e$, will be called \emph{good} for the chosen generating sequence.

\medskip

The following single-degree consequence is the form that will be used later in the apolar argument.

\begin{corollary}\label{thm:frob-amp}
Let $g_1,\ldots,g_s\in S$ be homogeneous and suppose $(g_1,\ldots,g_s)_A=S_A$ for some $A\ge0$. Then, outside finitely many rational primes, for every $q=p^e$,
$$
(g_1^q,\ldots,g_s^q)_{qA+(N-1)(q-1)}
=
S_{qA+(N-1)(q-1)}.
$$
\end{corollary}

\begin{proof}
Set $J=(g_1,\ldots,g_s)$ and
$
D=qA+(N-1)(q-1).
$
Since $J_A=S_A$, one has $(S/J)_m=0$ for all $m\ge A$. Thus a nonzero summand on the right side of \eqref{eq:HF-dilation-coeff} would require $m\le A-1$. For such $m$,
$$
\begin{aligned}
D-qm
&=q(A-m)+(N-1)(q-1)\\
&\ge q+(N-1)(q-1)\\
&=N(q-1)+1.
\end{aligned}
$$
Hence $C_{N,q}(D-qm)=0$ for every possible $m$, and \eqref{eq:HF-dilation-coeff} gives
$$
\dim_{\kk}(S/J_q)_D=0.
$$
Therefore $(J_q)_D=S_D$.
\end{proof}

\section{Polynomial interpolation and a Frobenius proof of Demailly's conjecture}\label{sec:fatpoints}

We first apply Theorem \ref{thm:HF-dilation} to polynomial interpolation and, particularly, to Demailly's conjecture. 
Throughout this section, assume that $n\ge1$.

\subsection{Inverse systems}

Let
$$
R=\kk[x_0,\ldots,x_n],
\quad
S=\kk[y_0,\ldots,y_n].
$$
We let $R$ act on $S$ by differentiation, so that $x_i$ acts as $\partial/\partial y_i$. For each $d$, this gives the usual perfect apolar pairing between $R_d$ and $S_d$. If $V\subseteq R_d$ is a $\kk$-vector subspace, we write
$$
V^\perp=\{G\in S_d:f\circ G=0\text{ for every }f\in V\}.
$$
Since $\kk$ has characteristic zero, this polynomial model identifies degree by degree with the divided power dual.

For a point $P_i\in\PP_\kk^n$, let $L_i\in S_1$ denote the corresponding dual linear form. The Emsalem--Iarrobino inverse system formula \cite[Theorem~I]{EI} gives, for positive integers $u_1,\ldots,u_s$ and every degree $D$,
\begin{equation}\label{eq:EI}
\left(\bigcap_{i=1}^s I(P_i)^{u_i}\right)_D^\perp
=
\sum_{i=1}^s L_i^{D-u_i+1}S_{u_i-1},
\end{equation}
provided $u_i\le D$ for every $i$. If $u_i>D$ for some $i$, then $(I(P_i)^{u_i})_D=0$, and hence the degree-$D$ piece of the intersection is zero. Thus, whenever $u_i\le D$ for all $i$,
\begin{equation}\label{eq:apolar-vanishing}
\left(\bigcap_{i=1}^s I(P_i)^{u_i}\right)_D=0
\quad\Longleftrightarrow\quad
\sum_{i=1}^s L_i^{D-u_i+1}S_{u_i-1}=S_D.
\end{equation}

The form of \eqref{eq:apolar-vanishing} used below is especially simple. Suppose
$$
F=\bigcap_{i=1}^s I(P_i)^{u_i}
\quad\text{and}\quad
a=\alpha(F).
$$
Since $F_{a-1}=0$, if $u_i<a$ for every $i$ then
\begin{equation}\label{eq:source-spanning}
S_{a-1}
=
\sum_{i=1}^s L_i^{a-u_i}S_{u_i-1}.
\end{equation}
In other words, the homogeneous ideal
$
J=(L_1^{a-u_1},\ldots,L_s^{a-u_s})
$
contains every form of degree $a-1$. 

\subsection{Quantitative postulation}


Let $P_1,\ldots,P_s\in\PP_\kk^n$ be distinct points and let $L_i\in S_1$ be their dual linear forms. For an integer vector $\mathbf u=(u_1,\ldots,u_s)$, use the convention $I(P_i)^{u_i}=R$ if $u_i\le0$, and put
$$
I(\mathbf u)=\bigcap_{i=1}^sI(P_i)^{u_i},
\quad
h_{\mathbf u}(D)=\dim_{\kk} I(\mathbf u)_D,
$$
with $h_{\mathbf u}(D)=0$ for $D<0$.

Fix $A\ge\max_i u_i\ge1$ and set
$$
r_i=A-u_i+1,
\quad
J=(L_1^{r_1},\ldots,L_s^{r_s})\subseteq S.
$$
By Emsalem--Iarrobino \cite[Theorem~I]{EI}, for every integer $j$ with $A+j\ge0$,
\begin{equation}\label{eq:EI-all-shifts}
\dim_{\kk}(S/J)_{A+j}
=
h_{\mathbf u+j\mathbf1}(A+j).
\end{equation}
Indeed, if $u_i+j>0$, then $u_i\le A$ gives
$
u_i+j\le A+j,
$
so the Emsalem--Iarrobino formula applies to that multiplicity, and the corresponding exponent in the inverse system is
$$
(A+j)-(u_i+j)+1=A-u_i+1=r_i.
$$
If $u_i+j\le0$, then $r_i>A+j$, so $L_i^{r_i}$ has no contribution in degree $A+j$, exactly matching the convention that a nonpositive multiplicity imposes no condition. Applying Emsalem--Iarrobino to the indices with $u_i+j>0$ gives \eqref{eq:EI-all-shifts}.

\begin{theorem}\label{thm:fatpoint-HF}
For every good prime power $q$ associated to $J$ and every integer $t$ with $D=qA+t\ge0$, define
$$
v_i(q,t)=q(u_i-1)+t+1.
$$
Then
\begin{equation}\label{eq:fatpoint-HF}
h_{\mathbf v(q,t)}(qA+t)
\le
\sum_{j\in\ZZ}
C_{n+1,q}(t-qj)
\,h_{\mathbf u+j\mathbf1}(A+j).
\end{equation}
Only finitely many terms on the right are nonzero.
\end{theorem}

\begin{proof}
The $q$th power-generated ideal associated to $J$ is
$
J_q=(L_1^{qr_1},\ldots,L_s^{qr_s}).
$
For $D=qA+t$ and $v_i=v_i(q,t)$,
$$
\begin{aligned}
D-v_i+1
&=qA+t-\bigl(q(u_i-1)+t+1\bigr)+1\\
&=q(A-u_i+1) =qr_i.
\end{aligned}
$$
If $v_i>0$, then $qr_i\ge1$, so
$
v_i=D-qr_i+1\le D,
$
and the Emsalem--Iarrobino formula applies to that multiplicity. If $v_i\le0$, then
$$
qr_i=D-v_i+1>D,
$$
so $L_i^{qr_i}$ has no contribution in degree $D$, exactly as the factor $I(P_i)^{v_i}=R$ imposes no condition. Applying Emsalem--Iarrobino to the indices with $v_i>0$ therefore gives
$$
\dim_{\kk}(S/J_q)_D=h_{\mathbf v(q,t)}(D).
$$

Applying \eqref{eq:HF-dilation-coeff} to $J$ yields
$$
\begin{aligned}
h_{\mathbf v(q,t)}(qA+t)
&\le
\sum_{m\ge0}C_{n+1,q}(qA+t-qm)\dim_{\kk}(S/J)_m.
\end{aligned}
$$
Write $m=A+j$. Then
$
qA+t-qm=t-qj,
$
and \eqref{eq:EI-all-shifts} gives
$$
\dim_{\kk}(S/J)_{A+j}=h_{\mathbf u+j\mathbf1}(A+j).
$$
Here $m\ge0$ corresponds to $j\ge-A$. Since $h_{\mathbf u+j\mathbf1}(A+j)=0$ for $j<-A$ by convention, we may extend the sum to all $j\in\ZZ$. Therefore
$$
h_{\mathbf v(q,t)}(qA+t)
\le
\sum_{j\in\ZZ}
C_{n+1,q}(t-qj)
\,h_{\mathbf u+j\mathbf1}(A+j),
$$
which is \eqref{eq:fatpoint-HF}. The sum is finite because $C_{n+1,q}(a)=0$ unless $0\le a\le(n+1)(q-1)$.
\end{proof}

\subsection{Frobenius amplification and a new proof of Demailly's conjecture}

For initial-degree questions, only the vanishing case of Theorem \ref{thm:fatpoint-HF} is needed. The resulting amplification statement also allows the multiplicities to vary independently from point to point. In the uniform case it recovers the estimate underlying \cite[Theorem~3.8]{HaSivakumar}.

\begin{theorem}\label{thm:fatpoint-amp}
Fix positive integers $u_1,\ldots,u_s$ and put
$$
F=\bigcap_{i=1}^s I(P_i)^{u_i},
\quad
a=\alpha(F).
$$
Outside finitely many rational primes there is the following property: for every good prime power $q=p^e$ and every positive vector $\mathbf v$ satisfying
\begin{equation}\label{eq:target-mult}
v_i\ge qu_i+(n-1)(q-1)
\quad(1\le i\le s),
\end{equation}
one has
\begin{equation}\label{eq:fatpoint-amp}
\alpha(I(\mathbf v))
\ge
qa+(n-1)(q-1).
\end{equation}
\end{theorem}

\begin{proof}
For every $i$,
$$
F\subseteq I(P_i)^{u_i}
\quad\Longrightarrow\quad
a=\alpha(F)\ge\alpha(I(P_i)^{u_i})=u_i.
$$
If $u_i=a$ for some $i$, then
$
I(\mathbf v)\subseteq I(P_i)^{v_i} \text{ and }
v_i\ge qa+(n-1)(q-1),
$
and hence
$$
\alpha(I(\mathbf v))
\ge v_i
\ge qa+(n-1)(q-1).
$$
We may therefore assume $u_i\le a-1$ for all $i$. Set
$$
A=a-1,
\quad
t=n(q-1).
$$
Then
$$
v_i(q,t)
=q(u_i-1)+n(q-1)+1
=qu_i+(n-1)(q-1),
$$
and
$
qA+t
=qa+(n-1)(q-1)-1.
$

We apply Theorem \ref{thm:fatpoint-HF}. If $j\ge0$, then $F_A=0$ gives $(S/J)_A=0$, hence
$$
(S/J)_{A+j}=0
\quad\Longrightarrow\quad
h_{\mathbf u+j\mathbf1}(A+j)=0.
$$
If $j<0$, then
$
t-qj \ge n(q-1)+q = (n+1)(q-1)+1,
$
so
$
C_{n+1,q}(t-qj)=0.
$
Thus every summand in \eqref{eq:fatpoint-HF} vanishes, and
$$
h_{\mathbf v(q,t)}(qA+t)=0.
$$
Therefore
$$
\alpha(I(\mathbf v(q,t)))
\ge qA+t+1
=qa+(n-1)(q-1).
$$
Finally, if $\mathbf v\ge\mathbf v(q,t)$ componentwise, then
$
I(\mathbf v)\subseteq I(\mathbf v(q,t)),
$
so
$$
\alpha(I(\mathbf v))
\ge
\alpha(I(\mathbf v(q,t)))
\ge qa+(n-1)(q-1)
$$
as desired.
\end{proof}

\begin{remark}\label{rem:prime-order}
Once the generating sequence is fixed, there is a finite set of bad primes independent of $q$. Hence, after choosing any good prime $p$, all powers $q=p^e$ may be used in the asymptotic arguments below.
\end{remark}

\begin{corollary}\label{cor:points}
Let $Z\subset\PP_\kk^n$ be a nonempty finite set of distinct points and let $I=I(Z)$. Then, for every $m\ge1$,
$$
\widehat\alpha(I)
\ge
\frac{\alpha(I^{(m)})+n-1}{m+n-1}.
$$
For $m=1$, this is Chudnovsky's bound.
\end{corollary}

\begin{proof}
For a finite reduced set of points,
$$
I^{(t)}=\bigcap_{P_i\in Z}I(P_i)^t
$$
for every $t\ge1$. Fix $m\ge1$ and put $a=\alpha(I^{(m)})$. Choose a good prime $p$ for the source multiplicity vector $(m,\ldots,m)$, and let $q=p^e$. Theorem \ref{thm:fatpoint-amp}, applied with
$
v_q=qm+(n-1)(q-1),
$
gives
$$
\alpha(I^{(v_q)})
\ge
qa+(n-1)(q-1).
$$
Since $v_q\to\infty$ and the Waldschmidt limit exists,
$$
\widehat\alpha(I)
=
\lim_{e\to\infty}\frac{\alpha(I^{(v_q)})}{v_q}
\ge
\lim_{e\to\infty}
\frac{qa+(n-1)(q-1)}{qm+(n-1)(q-1)}
=
\frac{a+n-1}{m+n-1},
$$
as claimed.
\end{proof}

Corollary \ref{cor:points} is Demailly's conjecture for finite reduced point sets. H\`a--Sivakumar recently proved the conjecture in full generality \cite[Theorem~1.2]{HaSivakumar}; the argument above gives a different proof in which the vanishing is obtained from the coefficientwise Hilbert-function dilation through apolarity.

\subsection{Weighted and subadditive multiplicity families}

The same argument as that of Corollary \ref{cor:points} applies when the multiplicities at different points grow at different positive subadditive rates.

For a graded family $\mathcal F=\{F_t\}_{t\ge1}$ for which the limit exists, write
$$
\widehat\alpha(\mathcal F)=\lim_{t\to\infty}\frac{\alpha(F_t)}{t}.
$$
For the families in the next theorem, the existence of the limit is automatic from subadditivity.

\begin{theorem}\label{thm:main}
Let $P_1,\ldots,P_s\in\PP_\kk^n$ be distinct points, let $\mu_i:\NN\to\NN$ be positive subadditive sequences, and set $\mathcal F=\{F_t\}_{t\ge1}$, where
$$
F_t=\bigcap_{i=1}^s I(P_i)^{\mu_i(t)},
\quad
\lambda_i=\lim_{t\to\infty}\frac{\mu_i(t)}t>0.
$$
Then, for every $m\ge1$,
$$
\widehat\alpha(\mathcal F)
\ge
\min_i
\frac{\lambda_i\bigl(\alpha(F_m)+n-1\bigr)}{\mu_i(m)+n-1}.
$$
\end{theorem}

\begin{proof}
First observe that $\mathcal F$ is a graded family. Indeed, for $r,t\ge1$, we have
$$
F_rF_t \subseteq \bigcap_{i=1}^s I(P_i)^{\mu_i(r)+\mu_i(t)} \subseteq \bigcap_{i=1}^s I(P_i)^{\mu_i(r+t)}
=F_{r+t},
$$
where the second inclusion follows from the subadditivity of each $\mu_i$. Consequently $\{\alpha(F_t)\}_{t \ge 1}$ is subadditive, so Fekete's lemma shows that the defining limit for $\widehat\alpha(\mathcal F)$ exists.

Fix $m\ge1$ and set $u_i=\mu_i(m)$. By Fekete's lemma,
$
\displaystyle \lambda_i
=
\lim_{t\to\infty}\frac{\mu_i(t)}t
=
\inf_{t\ge1}\frac{\mu_i(t)}t,
$
so
$$
\mu_i(t)\ge\lambda_it
\quad(t\ge1).
$$
Choose a good prime $p$ for the source vector $\mathbf u = (u_1, \dots, u_s)$. For $q=p^e$, define
$$
M_q
=
\max_i
\left\lceil
\frac{q\mu_i(m)+(n-1)(q-1)}{\lambda_i}
\right\rceil.
$$
Then, for every $i$,
$$
\mu_i(M_q)
\ge\lambda_iM_q
\ge q\mu_i(m)+(n-1)(q-1).
$$
Applying Theorem \ref{thm:fatpoint-amp} with multiplicities $u_i=\mu_i(m)$ and $\mu_i(M_q)$ gives
$$
\alpha(F_{M_q})
\ge
q\alpha(F_m)+(n-1)(q-1).
$$
Hence
$$
\frac{\alpha(F_{M_q})}{M_q}
\ge
\frac{q\alpha(F_m)+(n-1)(q-1)}{M_q}.
$$
Moreover,
$$
\frac{M_q}{q} = \max_i \frac1q \left\lceil \frac{q\mu_i(m)+(n-1)(q-1)}{\lambda_i} \right\rceil \longrightarrow
\max_i \frac{\mu_i(m)+n-1}{\lambda_i}.
$$
Since $M_q\to\infty$,
$$
\lim_{e\to\infty}
\frac{\alpha(F_{M_q})}{M_q}
=
\widehat\alpha(\mathcal F).
$$
Taking $e\to\infty$ in the preceding inequality yields
$$
\widehat\alpha(\mathcal F) \ge
\frac{\alpha(F_m)+n-1}
{\displaystyle\max_i\frac{\mu_i(m)+n-1}{\lambda_i}} =
\min_i
\frac{\lambda_i\bigl(\alpha(F_m)+n-1\bigr)}
{\mu_i(m)+n-1},
$$
which is the asserted inequality.
\end{proof}

The inequality in Theorem \ref{thm:main} is sharp in the stated generality: for one point $P$ and $F_t=I(P)^{\mu(t)}$, one has $\alpha(F_t)=\mu(t)$ and $\widehat\alpha(\mathcal F)=\lambda$, while the right side of the theorem is identically $\lambda$.

A particularly useful specialization is linear growth with fixed weights. Let $c_1,\ldots,c_s$ be positive integers and set
$$
F_t=\bigcap_{i=1}^s I(P_i)^{c_it},
\quad
c_{\min}=\min_i c_i.
$$

\begin{corollary}\label{cor:weighted}
For every $m\ge1$,
\begin{equation}\label{eq:weighted}
\widehat\alpha(\mathcal F)
\ge
\frac{c_{\min}\bigl(\alpha(F_m)+n-1\bigr)}{c_{\min}m+n-1}.
\end{equation}
Moreover, there exists a prime $p$ such that, for every $q=p^e$, if
$$
M_q
=
qm+
\left\lceil\frac{(n-1)(q-1)}{c_{\min}}\right\rceil,
$$
then
$$
\alpha(F_{M_q})
\ge
q\alpha(F_m)+(n-1)(q-1).
$$
\end{corollary}

\begin{proof}
Apply Theorem \ref{thm:main} with $\mu_i(t)=c_it$ and $\lambda_i=c_i$. Since
$$
c\longmapsto
\frac{c\bigl(\alpha(F_m)+n-1\bigr)}{cm+n-1}
$$
is nondecreasing for $c>0$, its minimum over the $c_i$ occurs at $c_{\min}$, giving \eqref{eq:weighted}.

For the second assertion, observe that for every $i$,
$
c_iM_q
\ge qc_im+(n-1)(q-1),
$
so Theorem \ref{thm:fatpoint-amp} gives
$$
\alpha(F_{M_q})
\ge q\alpha(F_m)+(n-1)(q-1)
$$
as claimed.
\end{proof}

The point-set inequality also extends to reduced closed subschemes with the ambient dimension in the denominator. The reduction is finite dimensional: in the degree just below $\alpha(I(X)^{(m)})$, only finitely many closed points of $X$ are needed to force the same vanishing.

\begin{corollary}\label{cor:subschemes}
Let $\varnothing\ne X\subsetneq\PP_\kk^n$ be a reduced closed subscheme with homogeneous radical ideal $I=I(X)$. Then, for every $m\ge1$,
$$
\widehat\alpha(I)
\ge
\frac{\alpha(I^{(m)})+n-1}{m+n-1}.
$$
\end{corollary}

\begin{proof}
Set $a=\alpha(I(X)^{(m)})$ and $A=a-1$. Since $X$ is reduced and $\kk$ has characteristic zero, Zariski--Nagata \cite{Zariski,Nagata}, together with the fact that differential powers commute with arbitrary intersections, gives
\begin{equation}\label{eq:subscheme-point-intersection}
I(X)^{(t)}=\bigcap_{P\in X}I(P)^t\quad(t\ge1),
\end{equation}
where $P$ ranges over the closed points of $X$; here $I(P)^{(t)}=I(P)^t$ for a point prime. Hence
$$
0=(I(X)^{(m)})_A=\bigcap_{P\in X}(I(P)^m)_A.
$$
Since $R_A$ is finite dimensional, finitely many of these subspaces already have zero intersection. Thus there is a finite set $Z\subset X$ such that $(I(Z)^{(m)})_A=0$, and therefore $\alpha(I(Z)^{(m)})\ge a$.

Since $I(X)\subseteq I(Z)$, monotonicity of differential powers and Zariski--Nagata give $I(X)^{(t)}\subseteq I(Z)^{(t)}$ for every $t$. For $t=m$ this yields $\alpha(I(Z)^{(m)})\le a$, so equality holds. For arbitrary $t$ it gives $\widehat\alpha(I(X))\ge\widehat\alpha(I(Z))$. Corollary \ref{cor:points} applied to $Z$ now gives
$$
\widehat\alpha(I(X))
\ge
\frac{\alpha(I(X)^{(m)})+n-1}{m+n-1},
$$
as claimed.
\end{proof}

Theorem \ref{thm:main} also applies to genuinely nonlinear multiplicity functions; for example, one may take $\mu_i(t)=\lceil c_it\rceil+b_i(t)$, where $c_i>0$ and $b_i:\NN\to\ZZ_{\ge0}$ are bounded subadditive sequences.

\subsection{Eulerian bounds for volumes on blow-ups}\label{sec:volumes}

We now return to the full postulation inequality of Theorem \ref{thm:fatpoint-HF}, which leads to a volume bound for divisors on blown up varieties.

Let
$
\pi:X=\operatorname{Bl}_{P_1,\ldots,P_s}\PP^n\longrightarrow\PP^n
$
be the blow-up at the points, let $H=\pi^*\mathcal O_{\PP^n}(1)$, and let $E_1,\ldots,E_s$ be the exceptional divisors. For nonnegative integers $m_i$,
$$
h^0\!\left(X,\mathcal O_X\left(dH-\sum_i m_iE_i\right)\right)
=
h_{\mathbf m}(d).
$$
Thus the section counts of these divisors are precisely fat-point postulation values, and their volumes record the corresponding asymptotic growth. The relation between positivity on blow-ups, Seshadri constants, and Nagata--SHGH type questions is classical; see, for example, \cite{BauerEtAl,DKMSHGH}. Balla\"y \cite{Ballay} relates volume functions on blow-ups to Seshadri constants and Nagata's conjecture. Here the Frobenius postulation inequality gives an upper bound for the volume using only finitely many neighboring postulation values.

\begin{corollary}\label{cor:finite-postulation}
Let $\mathbf u=(u_1,\ldots,u_s)\in\ZZ^s$ and let $A\ge\max_i u_i\ge1$. Set
$$
J=(L_1^{A-u_1+1},\ldots,L_s^{A-u_s+1}).
$$
Assume $q\ge n$ is a good prime power associated to $J$. Then
\begin{equation}\label{eq:finite-postulation}
h_{q\mathbf u}(q(A+1)-1)
\le
\sum_{\ell=0}^{n-1}
C_{n+1,q}((\ell+1)q-1)
\,h_{\mathbf u-\ell\mathbf1}(A-\ell).
\end{equation}
\end{corollary}

\begin{proof}
Take $t=q-1$ in Theorem \ref{thm:fatpoint-HF}. Then
$
\mathbf v(q,q-1)=q\mathbf u,
$
and the coefficient indexed by $j$ is
$
C_{n+1,q}(q-1-qj).
$

If $j\ge1$, then $q-1-qj<0$, so this coefficient vanishes. Write $j=-\ell$ for $j\le0$. The support condition becomes
$$
0\le(\ell+1)q-1\le(n+1)(q-1).
$$
The upper inequality is equivalent to
$
\ell q\le n(q-1).
$
If it holds, then
$$
\ell\le n-\frac nq<n,
$$
so $\ell\le n-1$. Conversely, when $0\le\ell\le n-1$ and $q\ge n$,
$$
\ell q\le(n-1)q\le n(q-1).
$$
Hence the support condition is equivalent to $0\le\ell\le n-1$, and \eqref{eq:fatpoint-HF} becomes exactly \eqref{eq:finite-postulation}.
\end{proof}

To pass from Corollary \ref{cor:finite-postulation} to a volume, we need the leading asymptotics of its coefficients. These limits are Eulerian numbers. Let $\eul{n}{\ell}$ denote the Eulerian number counting permutations of $\{1,\ldots,n\}$ with $\ell$ descents.

\begin{lemma}\label{lem:Eulerian}
For $0\le\ell\le n-1$,
\begin{equation}\label{eq:Eulerian-limit}
\lim_{q\to\infty}
\frac{n!}{q^n}
C_{n+1,q}((\ell+1)q-1)
=
\eul{n}{\ell}.
\end{equation}
\end{lemma}

\begin{proof}
By inclusion--exclusion,
$$
C_{n+1,q}((\ell+1)q-1)
=
\sum_{j=0}^{\ell}
(-1)^j\binom{n+1}{j}
\binom{(\ell+1-j)q+n-1}{n}.
$$
For each fixed $j$,
$$
\frac{n!}{q^n}
\binom{(\ell+1-j)q+n-1}{n}
\longrightarrow
(\ell+1-j)^n.
$$
Therefore
$$
\lim_{q\to\infty}
\frac{n!}{q^n}
C_{n+1,q}((\ell+1)q-1) =
\sum_{j=0}^{\ell}
(-1)^j\binom{n+1}{j}(\ell+1-j)^n = \eul{n}{\ell},
$$
where the last equality is the standard explicit formula for the Eulerian number; the $j=\ell+1$ term is zero.
\end{proof}

\begin{theorem}\label{thm:Eulerian-volume}
Let $u_1,\ldots,u_s$ be positive integers, let $A\ge\max_i u_i$, and set
$$
D=(A+1)H-\sum_{i=1}^su_iE_i.
$$
Then
\begin{equation}\label{eq:Eulerian-volume}
\vol_X(D)
\le
\sum_{\ell=0}^{n-1}
\eul{n}{\ell}
\,h_{\mathbf u-\ell\mathbf1}(A-\ell).
\end{equation}
Equivalently,
$$
\vol_X(D)
\le
\sum_{\ell=0}^{n-1}
\eul{n}{\ell}
h^0\!\left(
\PP^n,
\left(\bigcap_i\mathcal I_{P_i}^{\max\{u_i-\ell,0\}}\right)(A-\ell)
\right).
$$
\end{theorem}

\begin{proof}
If $D$ is not big, then $\vol_X(D)=0$ and there is nothing to prove. Assume $D$ is big. Choose a good prime $p$ for the ideal generated by powers of linear forms in Corollary \ref{cor:finite-postulation}, and let $q=p^e$. For $e\gg0$, one has $q\ge n$, hence
\begin{equation}\label{eq:volume-prelimit}
h_{q\mathbf u}(q(A+1)-1)
\le
\sum_{\ell=0}^{n-1}
C_{n+1,q}((\ell+1)q-1)
\,h_{\mathbf u-\ell\mathbf1}(A-\ell).
\end{equation}

Write $\mathcal I_{P_i}$ for the ideal sheaf of $P_i$. Since the points are distinct,
$$
\pi_*\mathcal O_X\left(-\sum_i qu_iE_i\right)
=
\bigcap_i\mathcal I_{P_i}^{qu_i}.
$$
Each $I(P_i)^{qu_i}$ is saturated, hence so is their finite intersection; it sheafifies to the displayed intersection. Therefore
$$
\begin{aligned}
h_{q\mathbf u}(q(A+1)-1)
&=\dim_{\kk}\left(\bigcap_i I(P_i)^{qu_i}\right)_{q(A+1)-1}\\
&=h^0\!\left(\PP^n,
\left(\bigcap_i\mathcal I_{P_i}^{qu_i}\right)(q(A+1)-1)\right).
\end{aligned}
$$
By the projection formula, the last term equals
$
h^0(X,\mathcal O_X(qD-H)).
$

Choose a hyperplane in $\PP^n$ avoiding all $P_i$, and let $Y\in|H|$ be its inverse image. Then
$$
Y\cong\PP^{n-1},
\quad
E_i\cap Y=\varnothing,
\quad
qD|_Y\cong\mathcal O_{\PP^{n-1}}(q(A+1)).
$$
From
$
0\to\mathcal O_X(qD-H)
\to\mathcal O_X(qD)
\to\mathcal O_Y(qD)
\to 0
$
we obtain
$$
\begin{aligned}
0
&\le h^0(X,qD)-h^0(X,qD-H) \le h^0(Y,qD|_Y)\\
&=\binom{q(A+1)+n-1}{n-1} = O(q^{n-1}).
\end{aligned}
$$
After extending scalars, if necessary, to an uncountable algebraically closed extension of $\kk$, the dimensions of all spaces of sections and the volume are unchanged. Thus Lazarsfeld--Musta\c{t}\u{a} \cite[Theorem~A]{LazarsfeldMustata} gives
$$
\vol_X(D)
=
\lim_{m\to\infty}\frac{n!}{m^n}h^0(X,mD).
$$
Restricting to $m=q=p^e$ and using the $O(q^{n-1})$ estimate,
$$
\begin{aligned}
\vol_X(D)
&=
\lim_{e\to\infty}\frac{n!}{q^n}h^0(X,qD)\\
&=
\lim_{e\to\infty}\frac{n!}{q^n}h^0(X,qD-H)\\
&=
\lim_{e\to\infty}\frac{n!}{q^n}h_{q\mathbf u}(q(A+1)-1).
\end{aligned}
$$
Multiplying \eqref{eq:volume-prelimit} by $n!/q^n$ and passing to the limit gives
$$
\begin{aligned}
\vol_X(D)
&\le
\sum_{\ell=0}^{n-1}
\left(
\lim_{e\to\infty}\frac{n!}{q^n}
C_{n+1,q}((\ell+1)q-1)
\right)
 h_{\mathbf u-\ell\mathbf1}(A-\ell)\\
&=
\sum_{\ell=0}^{n-1}
\eul{n}{\ell}\,h_{\mathbf u-\ell\mathbf1}(A-\ell),
\end{aligned}
$$
by Lemma \ref{lem:Eulerian}.
\end{proof}

The right-hand side of Theorem \ref{thm:Eulerian-volume} involves only the $n$ neighboring postulation values
$$
h_{\mathbf u-\ell\mathbf1}(A-\ell),
\quad 0\le\ell\le n-1.
$$
In particular, if all of these spaces vanish, then $\vol_X(D)=0$, so $D$ is not big.

\section{Powers of forms and complete ternary nondefectivity}\label{sec:syzygies}

We next turn to varieties of powers. By Terracini's lemma, their secant dimensions are controlled by tangent maps whose generators are powers of general forms. For the subgeneric secants it is enough to know that these generators have no syzygies up to one specified coefficient degree; determining their entire Hilbert series is unnecessary. The first result shows that Frobenius dilation enlarges precisely this initial syzygy range.

Suppose $g_1,\ldots,g_s\in S_r$ are nonzero forms of the same degree. Define
$$
\delta(g_1,\ldots,g_s)
=
\min\left\{t\ge0:\ker\left(S_t^s\longrightarrow S_{r+t}\right)\ne0\right\},
$$
where $(h_1,\ldots,h_s)$ is sent to $\sum_i h_ig_i$. If no such $t$ exists, set $\delta(g_1,\ldots,g_s)=\infty$.

\begin{theorem}\label{thm:syzygy-gap}
For every good prime power $q$ occurring in Theorem \ref{thm:HF-dilation},
\begin{equation}\label{eq:syzygy-gap}
\delta(g_1^q,\ldots,g_s^q)
\ge
q\,\delta(g_1,\ldots,g_s).
\end{equation}
Equivalently, if $0\le T<q\delta(g_1,\ldots,g_s)$, then
$$
S_T^s\longrightarrow S_{qr+T},
\quad
(h_i)_i\longmapsto\sum_i h_i g_i^q,
$$
is injective.
\end{theorem}

\begin{proof}
Write $\delta=\delta(g_1,\ldots,g_s)$. If $\delta=\infty$, then $s=1$, since for $s>1$ the Koszul relation $(g_2,-g_1,0,\ldots,0)$ has coefficient degree $r$. Since $S$ is a domain and $g_1^q\ne0$, multiplication by $g_1^q$ is injective in every degree, so $\delta(g_1^q)=\infty$. We may therefore assume $\delta<\infty$.

Set $J=(g_1,\ldots,g_s)$. For $0\le t<\delta$,
$$
\dim_{\kk}(S/J)_{r+t}=\dim_{\kk}S_{r+t}-s\dim_{\kk}S_t,
$$
so
$$
\Hilb_{S/J}(z)
\equiv
\frac{1-sz^r}{(1-z)^N}
\pmod{z^{r+\delta}}.
$$
After substituting $z^q$ and multiplying by $(1+z+\cdots+z^{q-1})^N=(1-z^q)^N/(1-z)^N$,
\begin{equation}\label{eq:truncated-target}
\Hilb_{S/J}(z^q)(1+z+\cdots+z^{q-1})^N
\equiv
\frac{1-sz^{qr}}{(1-z)^N}
\pmod{z^{q(r+\delta)}}.
\end{equation}
Fix $0\le T<q\delta$. Theorem \ref{thm:HF-dilation} and \eqref{eq:truncated-target} give
$$
\dim_{\kk}(S/J_q)_{qr+T}
\le
\dim_{\kk}S_{qr+T}-s\dim_{\kk}S_T.
$$
The reverse inequality is automatic because $(J_q)_{qr+T}$ is the image of $S_T^s\to S_{qr+T}$. Hence that map is injective, proving the theorem.
\end{proof}

We now specialize to three variables, where Anick's theorem supplies the generic Hilbert function needed to start the Frobenius argument. For $t\ge0$, set
$$
r_t=\dim_{\kk}\kk[x,y,z]_t=\binom{t+2}{2}.
$$

\subsection{Frobenius dilation and ternary secants}\label{subsec:secant-powers}

We shall now focus on the ternary Waring problem.

\begin{theorem}\label{thm:secant-tangent}
Let $q=p^e\ge2$, let $d\ge1$, and put $a=\lfloor d/q\rfloor$. If
\begin{equation}\label{eq:secant-source}
s r_a\le r_{d+a},
\end{equation}
then, for general $G_1,\ldots,G_s\in \kk[x,y,z]_d$, the map
$$
\kk[x,y,z]_d^s\longrightarrow \kk[x,y,z]_{(q+1)d},
\quad
(H_i)_i\longmapsto\sum_i H_iG_i^q
$$
is injective.
\end{theorem}

\begin{proof}
Work first over $K=\overline{\mathbb F}_p$. By Anick's theorem \cite{Anick}, general degree-$d$ forms $g_1,\ldots,g_s$ have Fr\"oberg's Hilbert function. Since
$$
\frac{r_{d+t}}{r_t}
=
\left(1+\frac{d}{t+1}\right)
\left(1+\frac{d}{t+2}\right),
$$
the ratio $r_{d+t}/r_t$ is strictly decreasing in $t$. Hence \eqref{eq:secant-source} gives $sr_t<r_{d+t}$ for $0\le t<a$ and $sr_a\le r_{d+a}$. Since $a<d$, only the linear term in the generators contributes to the raw Fr\"oberg coefficients in degrees $d+t$ with $0\le t\le a$. Thus the maps $K[x,y,z]_t^s\to K[x,y,z]_{d+t}$ are injective for $t<a$; if equality occurs at $t=a$, the source and target have the same dimension and the Fr\"oberg quotient is zero, so the map is an isomorphism there as well. Hence $\delta(g_1,\ldots,g_s)\ge a+1$.

Because $K$ is perfect, every homogeneous $H_i\in K[x,y,z]_T$ has a unique decomposition
$$
H_i=\sum_{\beta\in\{0,\ldots,q-1\}^3}h_{i,\beta}^{\,q}y^\beta.
$$
If $\sum_iH_ig_i^q=0$, freeness over the subring of $q$th powers gives
$$
\left(\sum_i h_{i,\beta}g_i\right)^q=0
$$
for every $\beta$. Frobenius is injective on the polynomial ring, so $\sum_i h_{i,\beta}g_i=0$. If $T<q(a+1)$, each nonzero $h_{i,\beta}$ has degree at most $a$, hence all vanish. Thus $\delta(g_1^q,\ldots,g_s^q)\ge q(a+1)>d$, and taking $T=d$ gives the desired injectivity in characteristic $p$.

The universal tangent matrix has entries in $\ZZ$-polynomials in the coefficients of the $G_i$. The characteristic-$p$ injectivity therefore gives a maximal minor whose reduction modulo $p$ is nonzero. That minor is a nonzero polynomial over $\ZZ$, so its nonvanishing defines a nonempty open set after base change to $\kk$. Hence the tangent map is injective for general $G_1,\ldots,G_s$ over $\kk$.
\end{proof}

Applying Theorem \ref{thm:secant-tangent} in the coefficient degree $d$ and then using Terracini's lemma gives the following subgeneric nondefectivity criterion.

\begin{corollary}\label{cor:secant-criterion}
Let $q=p^e\ge2$, put $a=\lfloor d/q\rfloor$, and assume
$$
s\binom{a+2}{2}\le\binom{d+a+2}{2},
\quad
s\binom{d+2}{2}\le\binom{(q+1)d+2}{2}.
$$
Then
$$
\dim\sigma_s(V^{q+1}_{2,d})
=
s\binom{d+2}{2}-1.
$$
\end{corollary}

The preceding criterion controls all subgeneric secants once the numerical inequality holds. In the ternary case, Fr\"oberg--Ottaviani--Shapiro \cite{FOS} provide a uniform filling bound at the next secant order, so the two ingredients combine to give complete nondefectivity in an explicit large-degree range.

\begin{corollary}\label{cor:frobenius-complete}
Let $k\ge3$ and assume that $q=k-1$ is a prime power. Put
$$
D_k=3k^3-6k^2+4.
$$
If $d\ge D_k$, then $V^k_{2,d}$ is completely nondefective.
\end{corollary}

\begin{proof}
Write $d=(k-1)a+t$ with $0\le t\le k-2$. By Corollary \ref{cor:secant-criterion}, it is enough to verify the source inequality for $s=k^2-1$. Set
$$
\Delta_k(a,t)
=
2\bigl(r_{ka+t}-(k^2-1)r_a\bigr).
$$
A direct expansion gives
$$
\Delta_k(a,t)
=
a^2-3k^2a+2kat+3ka+3a-2k^2+t^2+3t+4.
$$
Let
$$
A_k=3k^2-3k-3,
$$
so that $D_k=(k-1)A_k+1$. For fixed $a$, $\Delta_k(a,t)$ is increasing in $t$, and
$$
\Delta_k(A_k,1)=2(k-1)(k+1)(3k-4)>0,
\quad
\Delta_k(A_k+1,0)=(k-2)(k-1)>0.
$$
Moreover,
$$
\Delta_k(a+1,0)-\Delta_k(a,0)
=2a-3k^2+3k+4,
$$
which is positive for $a\ge A_k+1$. Because $D_k=(k-1)A_k+1$, the condition $d\ge D_k$ means either $a=A_k$ and $t\ge1$, or $a\ge A_k+1$. The preceding inequalities therefore give $\Delta_k(a,t)>0$ whenever $d\ge D_k$. Thus Corollary \ref{cor:secant-criterion} gives the expected dimension for every $s\le k^2-1$ once these secants are subgeneric. This is automatic here: since $D_k\ge3k^2$ for $k\ge3$,
$$
\begin{aligned}
2\bigl(r_{kd}-(k^2-1)r_d\bigr)
&=d\bigl(d-3(k^2-k-1)\bigr)-2k^2+4\\
&\ge 9k^3+7k^2+4>0.
\end{aligned}
$$
At the next order,
$$
k^2r_d-r_{kd}=\frac{(k-1)(3dk+2k+2)}2>0,
$$
so the $k^2$-nd secant is expected to fill. Fr\"oberg--Ottaviani--Shapiro \cite[Theorem~4]{FOS} prove over $\mathbb C$ that
$$
\sigma_{k^2}(V^k_{2,d})=\PP(\mathbb C[x,y,z]_{kd}).
$$
The secant parameterization is defined over $\mathbb Q$, and filling is equivalent to the nonvanishing of a maximal-rank minor of its differential. Such a minor nonzero over $\mathbb C$ remains a nonzero polynomial after base change to $\kk$. Thus $\sigma_{k^2}(V^k_{2,d})$ fills over $\kk$ as well. Hence the $k^2$-nd and all higher secants have the expected dimension, proving complete nondefectivity.
\end{proof}

\begin{remark}\label{rem:frobenius-thresholds}
Corollary \ref{cor:frobenius-complete} gives the two thresholds used below:
$$
k=3:\ d\ge31,
\quad
k=4:\ d\ge100.
$$
For fourth powers, writing $d=3m+r$, $r\in\{0,1,2\}$, the critical $s=15$ source inequality becomes
$$
m^2-33m-28\ge0,
\quad
m^2-25m-24\ge0,
\quad
m^2-17m-18\ge0,
$$
respectively. The uniform threshold in Corollary \ref{cor:frobenius-complete} is sharp for this direct criterion: at $d=D_k-1=(k-1)A_k$ one has
$$
\Delta_k(A_k,0)=-2(k^2-2)<0.
$$
\end{remark}

\subsection{Ternary cubes}\label{subsec:ternary-cubes}

For a form $F$ of degree $3d$ in three variables, its cube rank is the least $s$ for which
$$
F=G_1^3+\cdots+G_s^3
$$
with $G_i$ of degree $d$. The generic cube rank is therefore the least $s$ for which the secant variety $\sigma_s(V^3_{2,d})$ fills $\PP(\kk[x,y,z]_{3d})$. The generic $k$-rank conjecture of \cite[Conjecture~1.2]{LORS}, suggested by Ottaviani, predicts in the present case that
$$
\operatorname{rk}^{\circ}_3(3,3d)
=
\left\lceil\frac{\binom{3d+2}{2}}{\binom{d+2}{2}}\right\rceil.
$$
We prove the stronger statement that every secant variety of $V^3_{2,d}$ has the expected dimension. Corollary \ref{cor:frobenius-complete} already gives the range $d\ge31$. For the remaining degrees we specialize the tangent generators to powers of linear forms, so Emsalem--Iarrobino reduces the problem to the following calculation on del Pezzo surfaces. We state the calculation for all $d$, since that is its natural form.

\begin{lemma}\label{lem:cube-fat-points}
Let $1\le s\le8$, let $P_1,\ldots,P_s\in\PP^2$ be general points, and let
$
X_s=\operatorname{Bl}_{P_1,\ldots,P_s}\PP^2
$
be the blowup of $\mathbb{P}^2$ at the points. Write $H$ for the pullback of a line and $E_i$ for the exceptional curves. For $d\ge1$, set
$$
D_{s,d}=3dH-(d+1)\sum_{i=1}^sE_i.
$$
Then
$$
h^0(X_s,D_{s,d})
=
\max\{0,r_{3d}-sr_d\},
$$
except for
$
(s,d)=(5,2),(5,3),
$
where the dimensions are $1$ and $6$, respectively.
\end{lemma}

\begin{proof}
The case $s=1$ is immediate. Assume $2\le s\le8$. For general points, $X_s$ is a del Pezzo surface and its cone of curves is generated by the $(-1)$-curves. By the classification in \cite[Lemma~3]{FHH} (and restriction to the first $s$ exceptional classes), every nonexceptional $(-1)$-curve has the form
$$
C=aH-\sum b_iE_i,
\quad
\sum b_i=3a-1,
$$
with
$$
a\le1,2,3,6
$$
for $s\le4$, $s=5,6$, $s=7$, $s=8$, respectively. Therefore
$$
D_{s,d}\cdot C=d+1-3a,
\quad
D_{s,d}\cdot E_i=d+1,
$$
and $D_{s,d}$ is nef for
$$
d\ge2\ (s\le4),
\quad d\ge5\ (s=5,6),
\quad d\ge8\ (s=7),
\quad d\ge17\ (s=8).
$$
Since $-K_{X_s}$ is ample, $D_{s,d}-K_{X_s}=D_{s,d}+(-K_{X_s})$ is ample in these ranges. Hence Kodaira vanishing and Riemann--Roch give
$$
h^0(X_s,D_{s,d})=\chi(X_s,D_{s,d})=r_{3d}-sr_d.
$$

For the remaining cases we use the nef divisors
$$
F_4=2H-\sum_{1}^{4}E_i,
\quad
F_5=F_6=5H-2\sum E_i,
\quad
F_7=8H-3\sum E_i,
\quad
F_8=17H-6\sum E_i.
$$
Indeed, for the same curve $C$ their intersections are $1-a,2-a,3-a,6-a$, respectively, and their intersections with the exceptional curves are positive; the degree bounds above therefore prove nefness.

For $s=2,d=1$, the line $H-E_1-E_2$ is fixed and the residual class $2H-E_1-E_2$ has four sections. For $s=3,d=1$, the three lines through pairs of points are fixed and sum to $D_{3,1}$, so $h^0=1$. For $s=4,d=1$, one has $D_{4,1}\cdot F_4=-2$, so $D_{4,1}$ is not effective.

For $s=5$, let $C=2H-\sum_{1}^{5}E_i$, the unique conic through the five points. Then
$$
D_{5,2}=3C,
\quad
D_{5,3}=2C+F_5,
\quad
D_{5,4}=C+2F_5.
$$
The displayed copies of $C$ are fixed because the successive intersections with $C$ are negative. Thus $|D_{5,2}|=\{3C\}$ and $h^0(D_{5,2})=1$. The residual divisors $F_5$ and $2F_5$ are nef, so $F_5-K_{X_5}$ and $2F_5-K_{X_5}$ are ample; Kodaira vanishing and Riemann--Roch give
$$
h^0(X_5,F_5)=6,
\quad
h^0(X_5,2F_5)=16.
$$
Hence $h^0(D_{5,3})=6$ and $h^0(D_{5,4})=16$. The cases $d=2,3$ are precisely the two departures from the expected dimension, whereas $d=4$ has the expected value. Also $D_{5,1}\cdot F_5<0$, so $D_{5,1}$ is not effective.

For $s=6$, $D_{6,d}\cdot F_6=3d-12$, so $D_{6,d}$ is not effective for $d\le3$. At $d=4$, if
$$
Q_i=2H-\sum_{j\ne i}E_j,
$$
then the six classes $Q_i$ are effective $(-1)$-conics by the above classification; they are pairwise disjoint, satisfy $D_{6,4}=\sum_iQ_i$, and have $D_{6,4}\cdot Q_i=-1$. Each $Q_i$ is therefore a fixed component, and after removing all six nothing remains. Hence $h^0(D_{6,4})=1$.

For $s=7$, $D_{7,d}\cdot F_7=3d-21$, so $D_{7,d}$ is not effective for $d\le6$. At $d=7$, the seven classes
$$
C_i=3H-2E_i-\sum_{j\ne i}E_j
$$
are effective $(-1)$-curves by the classification, are pairwise disjoint, sum to $D_{7,7}$, and each has intersection $-1$ with it. Hence every $C_i$ is fixed and their sum exhausts the divisor class, so $h^0(D_{7,7})=1$.

For $s=8$, $D_{8,d}\cdot F_8=3d-48$, so $D_{8,d}$ is not effective for $d\le15$. At $d=16$, the eight classes
$$
C_i=6H-3E_i-2\sum_{j\ne i}E_j
$$
are effective $(-1)$-curves by the classification, are pairwise disjoint, sum to $D_{8,16}$, and each has intersection $-1$ with it. Hence every $C_i$ is fixed and their sum exhausts the divisor class, so $h^0(D_{8,16})=1$. These values agree with $\max\{0,r_{3d}-sr_d\}$ in every case not explicitly excepted.
\end{proof}

We now return from the fat-point calculation to the tangent maps of $V^3_{2,d}$. The lemma gives maximal rank after the specialization $G_i=L_i^d$ except for two pairs $(s,d)$; those two cases can be handled separately.

\begin{theorem}\label{thm:ternary-cubes-all}
For every $d\ge1$ and every $s\ge1$,
$$
\dim\sigma_s(V^3_{2,d})
=
\min\left\{
 s\binom{d+2}{2}-1,
 \binom{3d+2}{2}-1
\right\}.
$$
\end{theorem}

\begin{proof}
Let $T=\kk[x,y,z]$. The affine tangent space to the cone over $V^3_{2,d}$ at $G^3$ is $G^2T_d$. By Terracini's lemma, it is enough to show that for general $G_1,\ldots,G_s\in T_d$ the map
\begin{equation}\label{eq:cube-tangent-map}
\Phi_{s,d}:T_d^s\longrightarrow T_{3d},
\quad
(H_i)_i\longmapsto\sum_iH_iG_i^2
\end{equation}
has maximal rank.

If $d\ge31$, complete nondefectivity follows from Corollary \ref{cor:frobenius-complete} with $k=3$ and $q=2$. We may therefore assume $1\le d\le30$.

By \cite[Theorem~4]{FOS}, the ninth secant fills the ambient space over $\mathbb C$. The relevant secant parameterization is defined over $\mathbb Q$, so the filling statement transfers to $\kk$. It therefore suffices to treat $s\le8$. Specialize $G_i=L_i^d$ for general linear forms $L_i$, with dual points $P_i\in\PP^2$. The Emsalem--Iarrobino inverse-system formula \cite[Theorem~I]{EI} gives
$$
\dim_{\kk}\left(T/(L_1^{2d},\ldots,L_s^{2d})\right)_{3d}
=
h^0(X_s,D_{s,d}),
$$
because $3d-(d+1)+1=2d$. Outside the two exceptional pairs in Lemma \ref{lem:cube-fat-points}, the cokernel therefore has dimension $\max\{0,r_{3d}-sr_d\}$, exactly the cokernel dimension of a maximal-rank map $T_d^s\to T_{3d}$. Hence the specialization, and therefore a general tuple, has maximal rank except possibly for $(s,d)=(5,2),(5,3)$.

For $(s,d)=(5,3)$, maximal rank means injectivity of $T_3^5\to T_9$. Work over $K=\overline{\mathbb F}_2$, set $T_K=K[x,y,z]$, and let $g_1,\ldots,g_5$ be general cubics in $T_K$. Put $J=(g_1,\ldots,g_5)$. By Anick's theorem \cite{Anick},
$$
\dim_K(T_K/J)_3=5,
\quad
\dim_K(T_K/J)_4=0.
$$
The exact Frobenius formula in degree $9$ gives
$$
\dim_K(T_K/J^{[2]})_9
=
C_{3,2}(3)\cdot5+C_{3,2}(1)\cdot0
=5.
$$
Thus $\dim_K(J^{[2]})_9=50=5\dim_K(T_K)_3$, so $(T_K)_3^5\to (T_K)_9$ is injective. A nonzero maximal minor of the universal matrix is therefore nonzero modulo $2$, hence defines a nonzero polynomial over characteristic zero. Its nonvanishing is a nonempty open condition, proving maximal rank for general five-tuples of cubics over $\kk$.

For $(s,d)=(5,2)$, maximal rank means surjectivity of $T_2^5\to T_6$. Take
$$
G_1=xz+yz,
\quad G_2=y^2,
\quad G_3=xy+z^2,
\quad G_4=y^2+yz,
\quad G_5=x^2+xz.
$$
With the degree-two basis $z^2,yz,y^2,xz,xy,x^2$ and the degree-six monomials ordered first by increasing $x$-degree and then by increasing $y$-degree, form the $28\times30$ matrix of \eqref{eq:cube-tangent-map}. The $28\times28$ submatrix consisting of all six columns in the $G_1,G_2,G_3,G_5$ blocks and the four columns $z^2,yz,xz,x^2$ in the $G_4$ block has determinant $-16$. Since the determinant is nonzero, the displayed map is surjective. Surjectivity is open, hence holds for a general five-tuple.

This proves maximal rank for every $s\le8$ in the remaining range $d\le30$. The ninth secant, and therefore every higher secant, fills the ambient space. Together with the Frobenius range $d\ge31$, the claimed formula follows for every $d\ge1$.
\end{proof}

\begin{corollary}\label{cor:generic-ternary-cube-rank}
The generic cube rank of a ternary form of degree $3d$ is
$$
\operatorname{rk}^{\circ}_3(3,3d)
=
\left\lceil
\frac{\binom{3d+2}{2}}{\binom{d+2}{2}}
\right\rceil.
$$
Thus the $k=3$, three-variable case of the generic $k$-rank conjecture of \cite[Conjecture~1.2]{LORS}, suggested by Ottaviani, holds.
\end{corollary}

\begin{remark}
Theorem \ref{thm:ternary-cubes-all} is stronger than the generic-rank assertion because it determines every secant dimension. For comparison, the general degeneration theorem of Casarotti--Postinghel \cite{CasarottiPostinghel} gives no positive secant range when $k=3$.
\end{remark}

\subsection{Complete nondefectivity in three variables}\label{subsec:all-k-rank}

To pass from cubes to arbitrary $k\ge4$, we combine the Frobenius range with two plane-interpolation inputs. Recall that, in the ternary case,
$$
V^k_{2,d}
=
\{[G^k]:G\in\kk[x,y,z]_d\}
\subseteq
\PP(\kk[x,y,z]_{kd}).
$$
The first input is the following homogeneous postulation statement, which is sufficient for all of the low-degree cases needed below.

\begin{lemma}\label{lem:homogeneous-postulation}
Let $2\le m\le42$, let $P_1,\ldots,P_s\in\PP^2$ be general points, and let $D\ge3m$. Then
$$
\dim_{\kk}\left(\bigcap_{i=1}^s I(P_i)^m\right)_D
=
\max\{0,r_D-sr_{m-1}\}.
$$
Equivalently, the homogeneous system of degree-$D$ plane curves with multiplicity at least $m$ at $s$ general points is non-special.
\end{lemma}

\begin{proof}
We first work over $\mathbb C$. Ciliberto--Miranda classify the homogeneous $(-1)$-special systems \cite[Theorem~2.4]{CilibertoMiranda}: they occur only for $s=2,3,5,6,7,8$, and in each of the six listed families the degree is strictly smaller than $3m$. Indeed, the largest of their upper bounds is $(17m-2)/6<3m$. Hence a homogeneous system of degree $D\ge3m$ is not $(-1)$-special. Dumnicki's theorem \cite[Theorem~32]{Dumnicki} proves the homogeneous Harbourne--Hirschowitz conjecture for multiplicities $m\le42$, so in our range the system is non-special. Therefore
$$
\dim_{\mathbb C}\left(\bigcap_{i=1}^s I(P_i)^m\right)_D
=
\max\{0,r_D-sr_{m-1}\}.
$$
Finally, the assertion transfers to the algebraically closed characteristic-zero field $\kk$. The postulation condition is maximal rank of an interpolation matrix with integer-polynomial entries in the point coordinates. A maximal minor nonzero over $\mathbb C$ is a nonzero polynomial over $\mathbb Q$, hence remains nonzero after base change to $\kk$.
\end{proof}

For fourth powers, Lemma \ref{lem:homogeneous-postulation} settles the range $3\le d\le41$, while Corollary \ref{cor:frobenius-complete} settles every $d\ge100$. The second plane-geometric input bridges the intervening degrees by means of a uniform multipoint Seshadri estimate. For $s\ge1$, let
$$
\epsilon_s
=
\epsilon(\PP^2;P_1,\ldots,P_s)
$$
denote the multipoint Seshadri constant at $s$ very general points; equivalently, on $X_s=\operatorname{Bl}_{P_1,\ldots,P_s}\PP^2$,
$$
\epsilon_s
=
\sup\left\{c\ge0:H-c\sum_{i=1}^sE_i\text{ is nef}\right\}.
$$
For $0 < c < \epsilon_s$, the divisor $H-c\sum E_i$ is ample.

\begin{lemma}\label{lem:fourth-seshadri}
Let $d\ge39$ and $1\le s\le15$. For general $G_1,\ldots,G_s\in T_d$, where $T=\kk[x,y,z]$, the map
$$
T_d^s\longrightarrow T_{4d},
\quad
(H_i)_i\longmapsto\sum_iH_iG_i^3
$$
is injective.
\end{lemma}

\begin{proof}
We first work over $\mathbb C$. Specialize $G_i=L_i^d$ for very general linear forms $L_i$, with dual points $P_i\in\PP^2$, and write $X_s=\operatorname{Bl}_{P_1,\ldots,P_s}\PP^2$. By the Emsalem--Iarrobino formula, the cokernel dimension of the specialized map is
$$
h^0(X_s,D_{s,d}),
\quad
D_{s,d}=4dH-(d+1)\sum_{i=1}^sE_i.
$$
Since
$$
K_{X_s}=-3H+\sum_{i=1}^sE_i,
$$
we have
$$
D_{s,d}-K_{X_s}
=(4d+3)H-(d+2)\sum_{i=1}^sE_i.
$$
Harbourne--Ro\'e \cite[Remark~1.2]{HarbourneRoeSeshadri}, using Dumnicki's verification of the homogeneous Harbourne--Hirschowitz conjecture through multiplicity $42$ (see \cite[Theorem~32]{Dumnicki}), obtain

$$
\epsilon_{15}
>
\frac1{\sqrt{15}}\sqrt{1-\frac1{42\cdot15-84}}
=
\sqrt{\frac{109}{1638}}.
$$
The multipoint Seshadri constants are nonincreasing in the number of points, so $\epsilon_s\ge\epsilon_{15}$ for $s\le15$. Moreover, the function $(d+2)/(4d+3)$ is decreasing and
$$
\frac{109}{1638}-\left(\frac{41}{159}\right)^2
=
\frac{239}{4601142}>0.
$$
Hence, for every $d\ge39$ and $s\le15$,
$$
\frac{d+2}{4d+3}<\epsilon_s.
$$
It follows that $D_{s,d}-K_{X_s}$ is ample. Kodaira vanishing and Riemann--Roch give
$$
\begin{aligned}
h^0(X_s,D_{s,d})
&=\chi(X_s,D_{s,d})\\
&=r_{4d}-sr_d.
\end{aligned}
$$
The right side is positive for $d\ge39$ and $s\le15$, since
$$
r_{4d}-15r_d=\frac{d^2-33d-28}{2}>0.
$$
Thus the specialized tangent map has the expected cokernel dimension and is injective over $\mathbb C$. Equivalently, a maximal minor of the corresponding universal tangent matrix is a nonzero polynomial. The matrix has integer-polynomial entries in the coefficients of the $G_i$, so this minor is nonzero over $\mathbb Q$ and remains nonzero after base change to any algebraically closed field of characteristic zero. Its nonvanishing defines a nonempty open set on which the tangent map is injective. This proves the assertion over $\kk$.
\end{proof}

The cube theorem and the two interpolation lemmas above leave only the higher-power range in which the general degeneration theorem of Casarotti--Postinghel applies. Combining these inputs with the filling theorem of Fr\"oberg--Ottaviani--Shapiro gives the complete ternary statement.

\begin{theorem}\label{thm:ternary-k-rank}
For every $k\ge3$, every $d\ge2$, and every $s\ge1$,
\begin{equation}\label{eq:ternary-nondefectivity}
\dim\sigma_s(V^k_{2,d})
=
\min\left\{
 s\binom{d+2}{2}-1,
 \binom{kd+2}{2}-1
\right\}.
\end{equation}
Consequently,
\begin{equation}\label{eq:ternary-k-rank}
\operatorname{rk}^{\circ}_k(3,kd)
=
\left\lceil
\frac{r_{kd}}{r_d}
\right\rceil
=
\left\lceil
\frac{\binom{kd+2}{2}}{\binom{d+2}{2}}
\right\rceil.
\end{equation}
Thus the generic $k$-rank conjecture of \cite[Conjecture~1.2]{LORS}, suggested by Ottaviani, holds in three variables for every $k\ge3$, and in fact every ternary variety of powers $V^k_{2,d}$ with $k\ge3$ and $d\ge2$ is completely nondefective.
\end{theorem}

\begin{proof}
The case $k=3$ is Theorem \ref{thm:ternary-cubes-all}, so assume $k\ge4$. For general $G_1,\ldots,G_s\in T_d$, where $T=\kk[x,y,z]$, Terracini's lemma reduces the $s$th secant problem to the rank of
\begin{equation}\label{eq:k-tangent-map}
T_d^s\longrightarrow T_{kd},
\quad
(H_i)_i\longmapsto\sum_i H_iG_i^{k-1}.
\end{equation}
Specialize $G_i=L_i^d$ for general linear forms $L_i$, with dual points $P_i\in\PP^2$. The Emsalem--Iarrobino formula gives
\begin{equation}\label{eq:k-fatpoint-specialization}
\dim_{\kk}\left(T/(L_1^{d(k-1)},\ldots,L_s^{d(k-1)})\right)_{kd}
=
\dim_{\kk}\left(\bigcap_{i=1}^sI(P_i)^{d+1}\right)_{kd},
\end{equation}
because
$$
kd-(d+1)+1=d(k-1).
$$

We first treat $k=4$. Corollary \ref{cor:frobenius-complete}, with $q=3$, gives complete nondefectivity for every $d\ge100$, so it remains to consider $2\le d\le99$.

At $d=2$, one has $r_8/r_2=45/6$, so $s\le7$ is subgeneric and the eighth secant is expected to fill. The case $s=1$ is immediate. For $2\le s\le7$, put
$$
D_s=8H-3\sum_{i=1}^sE_i
$$
on $X_s$. For a nonexceptional $(-1)$-curve $C=aH-\sum b_iE_i$, the classification recalled in Lemma \ref{lem:cube-fat-points} gives $a\le3$, and hence
$$
D_s\cdot C=3-a\ge0,
\quad
D_s\cdot E_i=3.
$$
Thus $D_s$ is nef, so $D_s-K_{X_s}$ is ample because $-K_{X_s}$ is ample. Kodaira vanishing and Riemann--Roch give
$$
h^0(X_s,D_s)=r_8-sr_2=45-6s,
$$
and the tangent map has maximal rank for $s\le7$. For $s=8$, the nef divisor
$$
F_8=17H-6\sum_{i=1}^8E_i
$$
from Lemma \ref{lem:cube-fat-points} satisfies $D_8\cdot F_8=-8$, so $D_8$ is not effective. Hence the eighth secant fills, and therefore all higher secants fill as well.

For $3\le d\le41$, the multiplicity in \eqref{eq:k-fatpoint-specialization} is $m=d+1\le42$, and $4d\ge3(d+1)$. Lemma \ref{lem:homogeneous-postulation} applies for every $s$, so every secant variety of $V^4_{2,d}$ has the expected dimension.

For the remaining range $42\le d\le99$, Lemma \ref{lem:fourth-seshadri} shows that all secants of order $s\le15$ are nondefective. Since
$$
r_{4d}-15r_d=\frac{d^2-33d-28}{2}>0,
\quad
16r_d-r_{4d}=18d+15>0,
$$
the first fifteen secants are subgeneric and the sixteenth secant is expected to fill. By \cite[Theorem~4]{FOS},
$$
\sigma_{16}(V^4_{2,d})=\PP(T_{4d})
$$
over $\mathbb C$, and hence over $\kk$ by the same base-change argument used above. Therefore every secant of $V^4_{2,d}$ has the expected dimension. Together with Corollary \ref{cor:frobenius-complete}, this proves the fourth-power case for all $d\ge2$.

We next assume $k\ge5$. For $2\le d\le8$, one has $m=d+1\le9$ and
$$
kd\ge5d\ge3(d+1).
$$
Lemma \ref{lem:homogeneous-postulation} applies to \eqref{eq:k-fatpoint-specialization} for every $s$, so every secant variety of $V^k_{2,d}$ has the expected dimension.

It remains, for $k\ge5$, to consider $d\ge9$. Set $N=r_d-1\ge54$. Casarotti--Postinghel \cite[Theorem~1.1]{CasarottiPostinghel} prove that $V^k_{2,d}$ is not $s$-defective whenever
\begin{equation}\label{eq:CP-bound}
s\le B_k(N):=\frac1{N+1}\binom{N+k-3}{N}.
\end{equation}
Their theorem is stated over $\mathbb C$. For each fixed $s$ in this range, nondefectivity is equivalent, by Terracini's lemma, to the nonvanishing of a minor of the expected rank of the tangent matrix, whose entries are polynomials with integer coefficients in the coefficients of the forms $G_i$. Thus a minor that is nonzero over $\mathbb C$ is a nonzero polynomial over $\mathbb Q$ and remains nonzero after base change to our algebraically closed characteristic-zero field $\kk$. Hence the same nondefectivity bound holds over $\kk$.

For $k=5$,
$$
B_5(N)=\frac{N+2}{2}\ge28>24=5^2-1.
$$
Moreover,
$$
\frac{B_{k+1}(N)}{B_k(N)}=\frac{N+k-2}{k-2}.
$$
If $B_k(N)\ge k^2-1$, then $B_{k+1}(N)\ge(k+1)^2-1$ follows from
$$
(k^2-1)(N+k-2)\ge(k-2)k(k+2),
$$
whose left side minus right side is
$$
(N-2)k^2+3k+2-N>0
$$
for $N\ge54$ and $k\ge5$. Hence
$$
B_k(N)\ge k^2-1
$$
for every $k\ge5$. Thus all secants of order at most $k^2-1$ are nondefective. On the other hand, \cite[Theorem~4]{FOS} gives
$$
\sigma_{k^2}(V^k_{2,d})=\PP(T_{kd}),
$$
and
$$
k^2r_d-r_{kd}
=
\frac{(k-1)(3dk+2k+2)}2>0,
$$
so the expected dimension is already the ambient dimension at order $k^2$. Every higher secant also fills. This proves \eqref{eq:ternary-nondefectivity} for every $k\ge3$ and $d\ge2$.
\end{proof}

\section{Powers of general forms and Nicklasson's conjecture}\label{sec:Nicklasson-general}

The secant arguments of the preceding section use only a short initial range of the Hilbert function of powers of general forms. Nicklasson's conjecture asks for the stronger statement that the entire Hilbert series agrees with that of general forms of the corresponding degree. We begin by comparing the Frobenius-dilation inequality with the generic Hilbert function; in three variables this leads to two reductions of the conjecture and to a bounded strip containing every possible failure for prime-power exponents.

\subsection{Comparison with generic forms and Nenashev's theorem}

Let $H^{\rm gen}_{N,s,d}(m)$ denote the Hilbert function in degree $m$ of a quotient by $s$ general degree-$d$ forms in $N$ variables. Fr\"oberg and L\"ofwall showed that this Hilbert function is constant on a dense open set \cite{FrobergLofwall}; by semicontinuity, its value in each degree is the smallest attained value.

For general $G_1,\ldots,G_s\in S_d$ and every good prime power $q$ associated to this generating sequence, Theorem \ref{thm:HF-dilation} gives
\begin{equation}\label{eq:generic-sandwich}
H^{\rm gen}_{N,s,dq}(D)
\le
\dim_{\kk}\bigl(S/(G_1^q,\ldots,G_s^q)\bigr)_D
\le
\sum_{m\ge0}C_{N,q}(D-qm)H^{\rm gen}_{N,s,d}(m).
\end{equation}
Whenever the two bounds agree, the powers have the generic Hilbert function.

To make the upper bound in \eqref{eq:generic-sandwich} explicit in an initial range, we use Nenashev's higher-degree extension \cite[Theorem~1]{Nenashev} of the Hochster--Laksov maximal-rank phenomenon \cite{HochsterLaksov}. Only two standard $\operatorname{GL}$-stable cases are needed here: the full space $S_r$ and the Veronese variety $\{L^r:L\in S_1\}$.

\begin{theorem}\label{thm:Nenashev-range}
Let $\mathcal D_r$ be either $S_r$ or the Veronese variety $\{L^r:L\in S_1\}$, and let $g_1,\ldots,g_s$ be general elements of $\mathcal D_r$. Fix $a\ge1$ and assume
\begin{equation}\label{eq:Nenashev-condition}
s
\le
\frac{\dim_{\kk} S_{r+a}}{\dim_{\kk} S_a}-\dim_{\kk} S_a.
\end{equation}
Then $\delta(g_1,\ldots,g_s)\ge a+1$. Consequently, for every good prime power $q$ and every $0\le T<q(a+1)$,
\begin{equation}\label{eq:Nenashev-range}
\dim_{\kk}\bigl(S/(g_1^q,\ldots,g_s^q)\bigr)_{qr+T}
=
\dim_{\kk}S_{qr+T}-s\dim_{\kk}S_T.
\end{equation}
In these degrees the quotient has the same Hilbert function as a quotient by $s$ general degree-$qr$ forms.
\end{theorem}

\begin{proof}
Under \eqref{eq:Nenashev-condition}, Nenashev's theorem gives
$$
\dim_{\kk}(g_1,\ldots,g_s)_{r+a}=s\dim_{\kk}S_a,
$$
so $S_a^s\to S_{r+a}$ is injective. If a relation occurs in a smaller coefficient degree $b<a$, multiplying it by any nonzero form in $S_{a-b}$ gives a relation in degree $a$; hence no such relation exists. Thus $\delta(g_1,\ldots,g_s)\ge a+1$, and Theorem \ref{thm:syzygy-gap} gives \eqref{eq:Nenashev-range}. The corresponding injectivity conditions for a general degree-$qr$ tuple are open, giving the last assertion.
\end{proof}

The choice $\mathcal D_r=S_r$ gives intervals toward Nicklasson's conjecture, while the Veronese choice gives corresponding intervals for powers of general linear forms. For $a=1$, sharper results are supplied by Hochster--Laksov in the full-space case and, in the Veronese case, by Alexander--Hirschowitz through the inverse-system interpretation explained by Iarrobino; see \cite{AH, HochsterLaksov,IarrobinoII}.

\subsection{Three variables and Nicklasson's conjecture}\label{subsec:Nicklasson}

Let $S=\kk[y_1,y_2,y_3]$. For a formal series $A(z)=\sum_{j\ge0}a_jz^j$ with $a_0=1$, write $[A(z)]_+$ for the Fr\"oberg truncation: its coefficient in degree $j$ is $a_j$ as long as all coefficients through degree $j$ are positive, and is zero from the first nonpositive coefficient onward. By Anick's theorem \cite{Anick}, Fr\"oberg's predicted series \cite{Froberg} is the Hilbert series of a quotient by $s$ general degree-$e$ forms:
$$
\left[\frac{(1-z^e)^s}{(1-z)^3}\right]_+.
$$

Nicklasson \cite{Nicklasson} conjectured that powers of general forms have the same Hilbert series.

\begin{conjecture}[Nicklasson]\label{conj:Nicklasson}
Let $G_1,\ldots,G_s\in S_d$ be general forms with $d\ge2$, and let $m\ge2$. Then
\begin{equation}\label{eq:Nicklasson-conj}
\Hilb_{S/(G_1^m,\ldots,G_s^m)}(z)
=
\left[\frac{(1-z^{dm})^s}{(1-z)^3}\right]_+.
\end{equation}
Equivalently, $(G_1^m,\ldots,G_s^m)$ has the same Hilbert series as an ideal generated by $s$ general forms of degree $dm$.
\end{conjecture}

Because Anick's theorem determines exactly where the generic quotient first vanishes, the full conjecture in three variables is equivalent to two adjacent rank conditions. For $e\ge1$ and $s\ge4$, put
\begin{equation}\label{eq:rho-critical}
\rho(s,e)=\min\{t\ge0:s r_t\ge r_{e+t}\}.
\end{equation}
Since
$$
4r_{e-1}-r_{2e-1}=e>0,
$$
one has $\rho(s,e)\le e-1$ for every $s\ge4$.

\begin{proposition}\label{prop:Nicklasson-two-maps}
Fix $s\ge4$, $d\ge2$, and $m\ge2$, and put $e=dm$ and $\rho=\rho(s,e)$. For general $G_1,\ldots,G_s\in S_d$, define
$$
\Phi_t:S_t^s\longrightarrow S_{e+t},
\quad
(H_i)_i\longmapsto\sum_iH_iG_i^m.
$$
Then Nicklasson's conjecture for $(s,d,m)$ is equivalent to the following two conditions:
\begin{enumerate}
\item if $\rho>0$, the map $\Phi_{\rho-1}$ is injective;
\item the map $\Phi_\rho$ is surjective.
\end{enumerate}
\end{proposition}

\begin{proof}
Because $\rho\le e-1$, the degree-$(e+t)$ coefficient of the untruncated Fr\"oberg series for $s$ degree-$e$ forms is simply
$$
r_{e+t}-sr_t
$$
for $0\le t\le\rho$. It is positive for $t<\rho$ and nonpositive for $t=\rho$. By Anick's theorem \cite{Anick}, the generic quotient by $s$ degree-$e$ forms therefore has Hilbert function
$$
r_{e+t}-sr_t\quad(0\le t<\rho),
$$
and vanishes in every degree at least $e+\rho$.

Suppose first that the two stated rank conditions hold. If a nonzero relation among $G_1^m,\ldots,G_s^m$ occurred in coefficient degree $t<\rho-1$, multiplying it by any nonzero form of degree $\rho-1-t$ would give a nonzero relation in coefficient degree $\rho-1$, since $S$ is a domain. Thus every $\Phi_t$ is injective for $t<\rho$, and the quotient has the generic Hilbert function through degree $e+\rho-1$. Surjectivity of $\Phi_\rho$ gives
$$
(G_1^m,\ldots,G_s^m)_{e+\rho}=S_{e+\rho}.
$$
Since $S$ is standard graded, the same equality holds in all larger degrees. Hence \eqref{eq:Nicklasson-conj} follows. The converse is immediate from the generic Hilbert function just described.
\end{proof}

The two-map reduction is particularly effective for four generators, where a classical Lefschetz theorem gives the full Hilbert series.

\begin{proposition}\label{prop:Nicklasson-four}
Nicklasson's conjecture holds in three variables whenever $s\le4$.
\end{proposition}

\begin{proof}
For $s\le3$, the general forms $G_1,\ldots,G_s$ form a regular sequence, and so do their powers; the asserted Hilbert series is then the complete-intersection Hilbert series.

Let $s=4$ and put $e=dm$. Specialize $G_i=L_i^d$ for four general linear forms. After a linear change of coordinates, we may take
$$
L_1=y_1,
\quad L_2=y_2,
\quad L_3=y_3.
$$
Set
$$
A=S/(y_1^e,y_2^e,y_3^e).
$$
By the strong Lefschetz property for monomial complete intersections in characteristic zero \cite{Stanley}, multiplication by $L_4^e$ has maximal rank in every degree for a general linear form $L_4$. Thus the exact sequence defined by multiplication by $L_4^e$ gives, degree by degree,
$$
\dim_\kk(A/L_4^eA)_t
=
\max\{\dim_\kk A_t-\dim_\kk A_{t-e},0\}.
$$
Since $\Hilb_A(z)=(1-z^e)^3/(1-z)^3$, this is precisely
$$
\Hilb_{S/(L_1^e,L_2^e,L_3^e,L_4^e)}(z)
=
\left[\frac{(1-z^e)^4}{(1-z)^3}\right]_+.
$$
Once the quotient first vanishes, standard gradedness forces all subsequent pieces to vanish.

By Anick's theorem \cite{Anick}, the right side is the Hilbert series of four general degree-$e$ forms. The Hilbert function of a general four-tuple $(G_1^m,\ldots,G_4^m)$ is no larger than that of the specialization $(L_1^e,\ldots,L_4^e)$, while the generic degree-$e$ Hilbert function is a lower bound for every four-tuple of degree-$e$ forms. The two bounds agree, proving the result.
\end{proof}

For more generators the conjecture remains open, but for a prime-power exponent Frobenius confines every possible failure to a bounded strip of degrees. Its width is $2q-2$, depending only on the exponent $q$, although the location of the strip depends on $s$ and $d$. Write
\begin{equation}\label{eq:source-Froberg}
F_{s,d}(z)=\frac{(1-z^d)^s}{(1-z)^3}=\sum_{j\ge0}f_jz^j,
\quad
\tau(s,d)=\min\{j:f_j\le0\}.
\end{equation}
For $s\ge4$, the number $\tau(s,d)$ is finite.

\begin{theorem}\label{thm:Nicklasson-three-vars}
Let $s\ge4$ and $d\ge2$, let $G_1,\ldots,G_s\in S_d$ be general, and let $q=p^e\ge2$ be a prime power. Then
\begin{equation}\label{eq:Nicklasson-outside-strip}
\dim_{\kk}\bigl(S/(G_1^q,\ldots,G_s^q)\bigr)_D
=
[z^D]\left[\frac{(1-z^{dq})^s}{(1-z)^3}\right]_+
\end{equation}
for every
$$
D<q\tau(s,d)
\quad\text{or}\quad
D\ge q\tau(s,d)+2q-2.
$$
Consequently, Nicklasson's conjecture for exponent $q$ can fail only in the $2q-2$ degrees
\begin{equation}\label{eq:Nicklasson-strip}
q\tau(s,d)\le D\le q\tau(s,d)+2q-3.
\end{equation}
\end{theorem}

\begin{proof}
Work first over $K=\overline{\mathbb F}_p$ and choose general degree-$d$ forms $g_1,\ldots,g_s\in K[y_1,y_2,y_3]$. By Anick's theorem \cite{Anick},
$$
\Hilb_{K[y]/(g_1,\ldots,g_s)}(z)
=
[F_{s,d}(z)]_+.
$$
The exact Frobenius decomposition gives
\begin{equation}\label{eq:Nicklasson-charp-series}
\Hilb_{K[y]/(g_1^q,\ldots,g_s^q)}(z)
=
\bigl([F_{s,d}(z)]_+\bigr)\big|_{z\mapsto z^q}
(1+z+\cdots+z^{q-1})^3.
\end{equation}

Suppose $D<q\tau(s,d)$. If $C_{3,q}(D-qj)\ne0$, then $j<\tau(s,d)$, so the truncation on the right side of \eqref{eq:Nicklasson-charp-series} is invisible in degree $D$. Hence its degree-$D$ coefficient is
$$
[z^D]F_{s,d}(z^q)(1+z+\cdots+z^{q-1})^3
=
[z^D]\frac{(1-z^{dq})^s}{(1-z)^3}.
$$
Moreover, for every $0\le E\le D$, taking $j_E=\lfloor E/q\rfloor$ gives $C_{3,q}(E-qj_E)>0$ and $f_{j_E}>0$, while every other contributing $j$ is also smaller than $\tau(s,d)$. Thus all raw target coefficients through degree $D$ are positive, so the target Fr\"oberg truncation has not yet occurred. We have therefore obtained the generic degree-$dq$ Hilbert function in every degree $D<q\tau(s,d)$.

At the other end, the starting quotient in \eqref{eq:Nicklasson-charp-series} vanishes from degree $\tau(s,d)$ onward. Since $(1+z+\cdots+z^{q-1})^3$ has degree $3(q-1)$, the right side of \eqref{eq:Nicklasson-charp-series} has degree at most
$$
q(\tau(s,d)-1)+3(q-1)
=
q\tau(s,d)+2q-3.
$$
Thus the characteristic-$p$ quotient is zero in degree
$$
D_0=q\tau(s,d)+2q-2.
$$

It remains to lift these rank statements to characteristic zero. The multiplication matrices for the family $(G_1^q,\ldots,G_s^q)$ have entries which are integer polynomials in the coefficients of the $G_i$. For every $D<q\tau(s,d)$, the characteristic-$p$ calculation above gives exactly the Fr\"oberg--Anick Hilbert value for $s$ general degree-$dq$ forms; equivalently, the corresponding multiplication matrix has the maximal rank attained by a general degree-$dq$ tuple. Hence in each of these finitely many degrees there is a minor of that rank which is nonzero modulo $p$. In degree $D_0$ there is likewise a maximal minor giving surjectivity. The product of these minors is a nonzero polynomial modulo $p$, hence a nonzero polynomial over $\mathbb Z$. Its nonvanishing defines a nonempty open set in characteristic zero on which all these rank conditions hold simultaneously.

For a general tuple $(G_1,\ldots,G_s)$ in that open set, the resulting quotient dimensions in the low degrees are therefore at most the Fr\"oberg--Anick values. On the other hand, the Hilbert function of $s$ general degree-$dq$ forms is the minimum attained among all degree-$dq$ tuples, so those same values are universal lower bounds. Equality follows. Surjectivity in degree $D_0$ implies surjectivity in every larger degree, and hence both the power quotient and the generic degree-$dq$ quotient vanish there. This proves \eqref{eq:Nicklasson-outside-strip}.
\end{proof}

\begin{corollary}\label{cor:Nicklasson-squares}
For squares of general forms in three variables, Nicklasson's conjecture can fail only in the two degrees
$$
2\tau(s,d)
\quad\text{and}\quad
2\tau(s,d)+1.
$$
\end{corollary}

\begin{remark}
Proposition \ref{prop:Nicklasson-two-maps} and Corollary \ref{cor:Nicklasson-squares} give two complementary reductions of the square-power problem. Proposition \ref{prop:Nicklasson-two-maps} reduces Nicklasson's conjecture to the last injective and first surjective maps for forms of degree $2d$, whereas Corollary \ref{cor:Nicklasson-squares} confines any possible failure to the two degrees $2\tau(s,d)$ and $2\tau(s,d)+1$. These two pairs of boundary degrees need not coincide.
\end{remark}

\section{Further algebraic consequences of Frobenius dilation}\label{sec:further}

We record two further numerical consequences of the coefficientwise inequality. The first controls multiplicity under fixed-generator powers; in the Artinian case the same estimate refines to an inequality that recovers Lech's inequality asymptotically.

\subsection{Multiplicity and Lech's inequality}

\begin{corollary}\label{cor:multiplicity}
Assume $J$ is proper and put $c=\operatorname{ht}J$. For every good prime power $q$,
\begin{equation}\label{eq:multiplicity}
e(S/J_q)\le q^c e(S/J).
\end{equation}
The exponent and constant are sharp: equality holds when $s=c$ and the chosen generating sequence $g_1,\ldots,g_c$ is a homogeneous regular sequence.
\end{corollary}

\begin{proof}
Since $\sqrt{J_q}=\sqrt J$, the two quotients have the same dimension, say $d=N-c$. If $d=0$, then $c=N$ and multiplicity equals length. Evaluating \eqref{eq:HF-dilation} at $z=1$ gives
$$
e(S/J_q)=\ell(S/J_q)\le q^N\ell(S/J)=q^ce(S/J),
$$
which proves \eqref{eq:multiplicity} in this case.

Assume now that $d>0$. Write
$$
\Hilb_{S/J}(z)=\frac{P(z)}{(1-z)^d}.
$$
Then
$$
\Hilb_{S/J}(z^q)(1+\cdots+z^{q-1})^N
=
\frac{P(z^q)(1-z^q)^c}{(1-z)^N},
$$
whose Hilbert polynomial has multiplicity $q^cP(1)=q^ce(S/J)$. Since the two quotients have the same positive dimension, the coefficientwise inequality gives \eqref{eq:multiplicity} by comparing the leading coefficients of their Hilbert polynomials.

If $s=c$ and the chosen generating sequence is a homogeneous regular sequence, then so is its sequence of $q$th powers. The multiplicities are $\prod_i\deg g_i$ and $q^c\prod_i\deg g_i$, respectively, so equality holds.
\end{proof}


When $J$ is Artinian, the coefficientwise inequality contains more information than the multiplicity estimate alone and leads directly to Lech's inequality. Let $\mathfrak m=(y_1,\ldots,y_N)$ and assume that $J=(g_1,\ldots,g_s)$ is $\mathfrak m$-primary.

\begin{corollary}\label{cor:Lech}
For every good prime power $q$,
$$
\Hilb_{S/J^q}(z)
\preceq
\Hilb_{S/J_q}(z)
\preceq
\Hilb_{S/J}(z^q)(1+z+\cdots+z^{q-1})^N.
$$
In particular,
\begin{equation}\label{eq:finite-Lech}
\ell(S/J^q)+\ell(J^q/J_q)
=\ell(S/J_q)
\le q^N\ell(S/J).
\end{equation}
Consequently,
\begin{equation}\label{eq:Lech-regular}
e(JS_{\mathfrak m},S_{\mathfrak m})
\le N!\,\ell(S/J).
\end{equation}
\end{corollary}

\begin{proof}
The inclusion $J_q\subseteq J^q$ gives an exact sequence
$$
0\longrightarrow J^q/J_q
\longrightarrow S/J_q
\longrightarrow S/J^q
\longrightarrow0.
$$
Together with Theorem \ref{thm:HF-dilation}, this gives the coefficientwise inequalities and, after evaluating at $z=1$, \eqref{eq:finite-Lech}. Fix a good prime $p$ and put $q=p^e$. Since localization at $\mathfrak m$ preserves these finite lengths,
$$
\begin{aligned}
e(JS_{\mathfrak m},S_{\mathfrak m})
&=N!\lim_{e\to\infty}
\frac{\ell(S/J^{p^e})}{p^{eN}}\\
&\le N!\,\ell(S/J),
\end{aligned}
$$
which is \eqref{eq:Lech-regular}. Since $S_{\mathfrak m}$ is regular, this is the classical Lech inequality \cite{Lech}.
\end{proof}

The same exact sequence decomposes the finite-level defect as
$$
q^N\ell(S/J)-\ell(S/J^q)
=
\bigl(q^N\ell(S/J)-\ell(S/J_q)\bigr)
+\ell(J^q/J_q).
$$
Thus $J^q/J_q$ records the contribution of mixed products that are absent from the fixed-generator ideal $J_q$. Related ideals generated by $q$th powers of a fixed generating sequence occur in Ma--Smirnov's refinements of Lech's inequality \cite[Proof of Theorem~17]{MaSmirnov}.

A simple specialization of Theorem \ref{thm:HF-dilation} worth recording is obtained for ideals generated by powers of linear forms. Let
$$
Q_{\mathbf a}=(L_1^{a_1},\ldots,L_s^{a_s}),
\quad
\mathbf a=(a_1,\ldots,a_s)\in\NN^s.
$$
Taking $g_i=L_i^{a_i}$ gives, for every good prime power $q$,
$$
\Hilb_{S/Q_{q\mathbf a}}(z)
\preceq
\Hilb_{S/Q_{\mathbf a}}(z^q)(1+z+\cdots+z^{q-1})^N,
$$
where
$
Q_{q\mathbf a}=(L_1^{qa_1},\ldots,L_s^{qa_s}).
$
In particular, Corollary \ref{cor:multiplicity} applies whenever its hypothesis holds; if $Q_{\mathbf a}$ is $\mathfrak m$-primary, Corollary \ref{cor:Lech} gives the corresponding length inequalities.


\begin{thebibliography}{99}

\bibitem{Anick}
D. J. Anick,
\emph{Thin algebras of embedding dimension three},
J. Algebra \textbf{100} (1986), 235--259.

\bibitem{AH}
J. Alexander and A. Hirschowitz,
\emph{Polynomial interpolation in several variables},
J. Algebraic Geom. \textbf{4} (1995), 201--222.


\bibitem{Ballay}
F. Balla\"y,
\emph{Volume functions on blow-ups and Seshadri constants},
Proc. Amer. Math. Soc. \textbf{150} (2022), no.~5, 1925--1935.

\bibitem{BauerEtAl}
T. Bauer, S. Di Rocco, B. Harbourne, M. Kapustka, A. L. Knutsen,
W. Syzdek, and T. Szemberg,
\emph{A primer on Seshadri constants},
in \emph{Interactions of classical and numerical algebraic geometry},
Contemp. Math. \textbf{496}, Amer. Math. Soc., Providence, RI, 2009, 33--70.

\bibitem{BCMO}
A. T. Blomenhofer, A. Casarotti, M. Micha{\l}ek, and A. Oneto,
\emph{Identifiability for mixtures of centered Gaussians and sums of powers of quadratics},
Bull. Lond. Math. Soc. \textbf{55} (2023), no.~5, 2407--2424.

\bibitem{CasarottiPostinghel}
A. Casarotti and E. Postinghel,
\emph{Waring identifiability for powers of forms via degenerations},
Proc. Lond. Math. Soc. (3) \textbf{128} (2024), no.~1, e12579.


\bibitem{CilibertoMiranda}
C. Ciliberto and R. Miranda,
\emph{Linear systems of plane curves with base points of equal multiplicity},
Trans. Amer. Math. Soc. \textbf{352} (2000), no.~9, 4037--4050.

\bibitem{Chudnovsky}
G. V. Chudnovsky,
\emph{Singular points on complex hypersurfaces and multidimensional Schwarz lemma},
in S\'eminaire de Th\'eorie des Nombres, Paris 1979--80,
Progr. Math. \textbf{12}, Birkh\"auser, 1981, 29--69.

\bibitem{Demailly}
J.-P. Demailly,
\emph{Formules de Jensen en plusieurs variables et applications arithm\'etiques},
Bull. Soc. Math. France \textbf{110} (1982), 75--102.

\bibitem{Dumnicki}
M. Dumnicki,
\emph{Cutting diagram method for systems of plane curves with base points},
Ann. Polon. Math. \textbf{90} (2007), no.~2, 131--143.



\bibitem{DKMSHGH}
M. Dumnicki, A. K\"uronya, C. Maclean, and T. Szemberg,
\emph{Rationality of Seshadri constants and the Segre--Harbourne--Gimigliano--Hirschowitz conjecture},
Adv. Math. \textbf{303} (2016), 1162--1170.

\bibitem{EI}
J. Emsalem and A. Iarrobino,
\emph{Inverse system of a symbolic power. I},
J. Algebra \textbf{174} (1995), 1080--1090.



\bibitem{FHH}
S. Fitchett, B. Harbourne, and S. Holay,
\emph{Resolutions of fat point ideals involving eight general points of $\PP^2$},
J. Algebra \textbf{244} (2001), no.~2, 684--705.

\bibitem{HarbourneRoeSeshadri}
B. Harbourne and J. Ro\'e,
\emph{Computing multi-point Seshadri constants on $\PP^2$},
Bull. Belg. Math. Soc. Simon Stevin \textbf{16} (2009), no.~5, 887--906.

\bibitem{Froberg}
R. Fr\"oberg,
\emph{An inequality for Hilbert series of graded algebras},
Math. Scand. \textbf{56} (1985), 117--144.

\bibitem{FrobergLofwall}
R. Fr\"oberg and C. L\"ofwall,
\emph{On Hilbert series for commutative and noncommutative graded algebras},
J. Pure Appl. Algebra \textbf{76} (1991), 33--38.

\bibitem{FOS}
R. Fr\"oberg, G. Ottaviani, and B. Shapiro,
\emph{On the Waring problem for polynomial rings},
Proc. Natl. Acad. Sci. USA \textbf{109} (2012), no.~15, 5600--5602.

\bibitem{HaSivakumar}
T. H. H\`a and A. Sivakumar,
\emph{Polynomial interpolation and the Waldschmidt constant of points in projective space},
preprint, arXiv:2608.16040v1 (2026).

\bibitem{HochsterLaksov}
M. Hochster and D. Laksov,
\emph{The linear syzygies of generic forms},
Comm. Algebra \textbf{15} (1987), 227--239.

\bibitem{IarrobinoII}
A. Iarrobino,
\emph{Inverse system of a symbolic power. II. The Waring problem for forms},
J. Algebra \textbf{174} (1995), 1091--1110.

\bibitem{LazarsfeldMustata}
R. Lazarsfeld and M. Musta\c{t}\u{a},
\emph{Convex bodies associated to linear series},
Ann. Sci. \'Ec. Norm. Sup\'er. (4) \textbf{42} (2009), 783--835.

\bibitem{Lech}
C. Lech,
\emph{Note on multiplicities of ideals},
Ark. Mat. \textbf{4} (1960), 63--86.


\bibitem{LORS}
S. Lundqvist, A. Oneto, B. Reznick, and B. Shapiro,
\emph{On generic and maximal $k$-ranks of binary forms},
J. Pure Appl. Algebra \textbf{223} (2019), no.~5, 2062--2079.

\bibitem{MaSmirnov}
L. Ma and I. Smirnov,
\emph{Uniform Lech's inequality},
Proc. Amer. Math. Soc. \textbf{151} (2023), 2387--2397.

\bibitem{Nagata}
M. Nagata,
\emph{Local Rings},
Interscience Tracts in Pure and Applied Mathematics, Vol.~13,
Interscience Publishers, New York--London, 1962.

\bibitem{Nenashev}
G. Nenashev,
\emph{A note on Fr\"oberg's conjecture for forms of equal degrees},
C. R. Math. Acad. Sci. Paris \textbf{355} (2017), 272--276.

\bibitem{Nicklasson}
L. Nicklasson,
\emph{On the Hilbert series of ideals generated by generic forms},
Comm. Algebra \textbf{45} (2017), 3390--3395.

\bibitem{Stanley}
R. P. Stanley,
\emph{Weyl groups, the hard Lefschetz theorem, and the Sperner property},
SIAM J. Algebraic Discrete Methods \textbf{1} (1980), no.~2, 168--184.

\bibitem{Zariski}
O. Zariski,
\emph{A fundamental lemma from the theory of holomorphic functions on an algebraic variety},
Ann. Mat. Pura Appl. (4) \textbf{29} (1949), 187--198.

\end{thebibliography}
\end{document}